\documentclass[12pt]{amsart}
\usepackage{amssymb}

\newtheorem{theorem}{Theorem}
\newtheorem{example}{Example}
\newtheorem{proposition}{Proposition}
\newtheorem{corollary}{Corollary}
\newtheorem{remark}{Remark}
\newtheorem{definition}{Definition}
\newtheorem{lemma}{Lemma}

\newcommand{\Z}{{\mathbb Z}}
\newcommand{\Q}{{\mathbb Q}}
\newcommand{\R}{{\mathbb R}}
\newcommand{\C}{{\mathbb C}}
\newcommand{\PP}{{\mathfrak P}}
\newcommand{\RP}{{\mathbb RP}}
\newcommand{\CP}{{\mathbb CP}}

\begin{document}

\author{V.A.~Vassiliev}
\address{Steklov Mathematical Institute of Russian Academy of Sciences \ \ and \ \ National Research University Higher School of Economics}
 \email{vva@mi-ras.ru}

\thanks{}

\title[Homology of spaces of equivariant maps]{Twisted homology of configuration spaces, homology of spaces of equivariant maps, and stable homology of spaces of non-resultant systems of real homogeneous polynomials}

\begin{abstract}
A spectral sequence calculating the homology groups of some spaces of maps equivariant under compact group actions is described. For the main example, we calculate the rational homology groups of spaces of even and odd maps $S^m \to S^M$, $m<M$, or, which is the same, the stable homology groups of spaces of non-resultant homogeneous polynomial maps $\R^{m+1} \to \R^{M+1}$ of growing degrees. Also, we find the homology groups of spaces of $\Z_r$-equivariant maps of odd-dimensional spheres for any $r$.
In an intermediate calculation, we find the homology groups of configuration spaces of projective and lens spaces with coefficients in several local systems. 
\end{abstract}



\maketitle

\section{Main results}

\subsection{Homology of spaces of $\Z_r$-equivariant maps of spheres}
Denote by $\mbox{EM}_0(S^m,S^M)$ (respectively, $\mbox{EM}_1(S^m,S^M)$) the space of even (respectively, odd) continuous maps $S^m \to S^M$, that is, of maps sending opposite points to the same point (respectively, to opposite points) in $S^M$. (The space $\mbox{EM}_0(S^m,S^M)$ can be identified with the space of continuous maps $\RP^m \to S^M$). Also, denote by $\mbox{EM}^*_i(S^m,S^M)$, $i \in \{0, 1\},$ the pointed versions of these spaces, i.e. the spaces of even or odd maps sending the fixed point of $S^m$ into the fixed point of $S^M$.

\begin{theorem}
\label{cor1} For any natural $m<M$ and $i =0$ $($respectively, $i=1)$, the Poincar\'e series of the
group $H^*(\mbox{EM}_i(S^m,S^M), \Q)$ is indicated in the intersection of the corresponding row and the third column of Table
\ref{evv} $($respectively, of Table \ref{oddd}$)$. The Poincar\'e series of the similar group $H^*(\mbox{EM}^*_i(S^m,S^M), \Q)$ is indicated in the intersection of the corresponding row and the second column of the same table.
\end{theorem}


There are obvious fiber bundles 
\begin{equation}
\label{mfb}
\mbox{EM}_i (S^m, S^M) \to S^M
\end{equation}
 sending any map to the image of the fixed point of $S^m$ under this map; their fibers are equal to $\mbox{EM}^*_i(S^m,S^M)$. The last columns of Tables \ref{evv} and \ref{oddd} complete the description of the rational cohomology spectral sequences of these fiber bundles. Namely, it follows from Theorem \ref{cor1} that all the cohomology groups $H^q(\mbox{EM}^*_i(S^m,S^M), \Q)$ of their fibers are 
 at most one-dimensional. Therefore, to find the shape of all terms ${\mathcal E}_r$, $r\geq 1$, of such a spectral sequence, it is enough to indicate all the pairs of its groups ${\mathcal E}_M^{0,q}$ and ${\mathcal E}_M^{M, q-M+1}$, which are connected by a non-trivial action of the differential $d^M$. These pairs are listed in the last columns of our tables.

In Table \ref{genn} we similarly describe the spectral sequences of the fiber bundles 
\begin{equation}\label{fb}\mbox{Map}(S^m,S^M) \stackrel{\Omega^m S^M}{\longrightarrow} S^M
\end{equation}
for the spaces of all continuous maps $S^m \to S^M$, $m<M$.
The answers given in the second column of this table (and probably also in entire this table) are not new since \cite{snaith}, \cite{Cohen76} and related works. We give them here for the comparison with the data of two previous tables; also, some homology calculations arising in this approach to these answers may be of independent interest.

\begin{table}
\begin{center}
\begin{tabular}{|l|c|c|c|}
\hline
Parities & $\mbox{EM}^*_0(S^m, S^M)$ & $\mbox{EM}_0(S^m,S^M)$ & $d^M$ \\
\hline
$\mbox{$\begin{array}{l} $m$ \ odd, \\$M$ \ odd \end{array}$}$ & {\large $\frac{1}{1-t^{M-m}}$} & {\large $\frac{1+t^M}{1-t^{M-m}}$} & $\emptyset$ \\
\hline
$\mbox{$\begin{array}{c} $m$ \ even, \\$M$ \ even \end{array}$}$ & {\large 1} & {\large $1+ t^M $} & 
$\emptyset$ \\
\hline
$\mbox{$\begin{array}{c} $m$ \ even, \\$M$ \ odd \end{array}$}$ & {\large 1} & {\large $1+t^M$ } & $\emptyset$ \\
\hline 
$\mbox{$\begin{array}{c} $m$ \ odd, \\$M$ \ even \end{array}$}$
 & {\large $\frac{1+t^{M-m}}{1-t^{2M-m-1}}$} & {\large $1+ t^{M-m} \frac{1+t^m}{1-t^{2M-m-1}}$} & ${\mathcal E}_M^{0, (s+1)(2M-m-1)} \to {\mathcal E}_M^{M, s(2M-m-1)+(M-m)}$, $s \geq 0$ \\
\hline
\end{tabular}
\caption{Poincar\'e series of spaces of even maps $S^m \to S^M$}
\label{evv}
\end{center}
\end{table}

\begin{table}
\begin{center}
\begin{tabular}{|l|c|c|c|}
\hline
Parities & $\mbox{EM}^*_1(S^m, S^M)$ & $\mbox{EM}_1(S^m,S^M)$ & $d^M$ \\
\hline
$\mbox{$\begin{array}{c} $m$ \ odd, \\$M$ \ odd \end{array}$}$ & {\large $\frac{1}{1-t^{M-m}}$} & {\large $\frac{1+t^M}{1-t^{M-m}}$} & $\emptyset$ \\
\hline
$\mbox{$\begin{array}{c} $m$ \ even, \\$M$ \ even \end{array}$}$ & {\large $\frac{1+t^{M-1}}{1-t^{M-m}}$} & {\large $\frac{1+t^{2M-1}}{1-t^{M-m}}$} & ${\mathcal E}_M^{0,s(M-m)+M-1} \to {\mathcal E}_M^{M,s(M-m)}$, $s \geq 0$
 \\
\hline
$\mbox{$\begin{array}{c} $m$ \ even, \\$M$ \ odd \end{array}$}$ & {\large 1} & {\large $1+t^M$} & $\emptyset$ \\
\hline 
$\mbox{$\begin{array}{c} $m$ \ odd, \\$M$ \ even \end{array}$}$
 & {\large $\frac{1+t^{M-1}}{1-t^{2M-m-1}}$} & {\large $ \frac{1+t^{2M-1}}{1-t^{2M-m-1}}$} & ${\mathcal E}_M^{0,s(2M-m-1)+M-1} \to {\mathcal E}_M^{M,s(2M-m-1)}$, $s \geq 0$ \\
\hline
\end{tabular}
\caption{Poincar\'e series of spaces of odd maps $S^m \to S^M$}
\label{oddd}
\end{center}
\end{table}

\begin{table}
\begin{center}
\begin{tabular}{|l|c|c|c|}
\hline
Parities & $\Omega^m S^M$ & $\mbox{Map}(S^m,S^M)$ & $d^M$ \\
\hline
$\mbox{$\begin{array}{c} $m$ \ odd, \\$M$ \ odd \end{array}$}$ & {\large $\frac{1}{1-t^{M-m}}$} & {\large $\frac{1+t^M}{1-t^{M-m}}$} & $\emptyset$ \\
\hline
$\mbox{$\begin{array}{c} $m$ \ even, \\$M$ \ even \end{array}$}$ & {\large $\frac{1+t^{2M-m-1}}{1-t^{M-m}}$} & {\large $t^M + \frac{1+t^{3M-m-1}}{1-t^{M-m}}$} & 
${\mathcal E}_M^{0,(s-1)(M-m)+(2M-m-1)} \to {\mathcal E}_M^{M, s(M-m)} ,$ $s \geq 1$ \\
\hline
$\mbox{$\begin{array}{c} $m$ \ even, \\$M$ \ odd \end{array}$}$ & $1+t^{M-m}$ & $(1+t^M)(1+t^{M-m})$ & $\emptyset$ \\
\hline 
$\mbox{$\begin{array}{c} $m$ \ odd, \\$M$ \ even \end{array}$}$
 & {\large $\frac{1+t^{M-m}}{1-t^{2M-m-1}}$} & {\large $1+ t^{M-m} \frac{1+t^m}{1-t^{2M-m-1}}$} & ${\mathcal E}_M^{0, (s+1)(2M-m-1)} \to {\mathcal E}_M^{M, s(2M-m-1)+(M-m)}$, $s \geq 0$ \\
\hline
\end{tabular}
\caption{Poincar\'e series of spaces of general maps $S^m \to S^M$}
\label{genn}
\end{center}
\end{table}

Any odd-dimensional sphere can be realized as the unit sphere in a complex vector space, hence the cyclic group $\Z_r$ acts on it by multiplications by the powers of $e^{2\pi i/r}$. 

\begin{theorem} \label{m4}
Let $m<M$ be two odd numbers, and $0 < s \leq r$, then 

1$)$ the Poincar\'e series of the rational cohomology group of the space of continuous maps $f: S^m \to S^M$ 
such that 
\begin{equation}
\label{chara}
f(e^{2\pi i/r} x) = e^{2\pi i s/r} f(x)
\end{equation}
 for any $x \in S^m$ is equal to $(1+t^{M})/(1-t^{M-m})$; 

2$)$ the Poincar\'e series of the corresponding space of equivariant maps sending the fixed point of $S^m$ to the fixed of $S^M$ is equal to $1/(1-t^{M-m})$.
\end{theorem} 

See \cite{circle} for the first application of our spectral sequence to the case of infinite group actions. \medskip

The data of Tables 1, 2 and 3 for $m$ and $M$ odd are special cases of Theorem 2 corresponding to $(r,s)$ equal to $(2,2),$ $(2,1)$, and $(1,1)$ respectively.

\subsection{Application to the spaces of non-resultant maps}
 
Denote by $\PP_{d,n}$ the space of homogeneous polynomials $\R^n \to \R$ of degree $d$. Its $k$-th Cartesian power
$\PP_{d,n}^k$ is the space of systems of $k$ such polynomials. The {\em resultant} subvariety 
$\Sigma \subset \PP_{d,n}^k$ consists of all systems having non-trivial common roots in $\R^n$.

The study of the topology of non-resultant spaces was essentially started in \cite{CCMM}; for some other related works see \cite{Book}, \cite{Kozl}, \cite{KY}, \cite{dan18}.
 
The embedding $\PP_{d,n}^k \to \PP_{d+2,n}^k$ defined by the multiplication of all polynomials by
$x_1^2 + \dots + x_n^2$ induces a homomorphism
\begin{equation} H^*(\PP^k_{d+2,n} \setminus \Sigma) \to
H^*(\PP^k_{d,n} \setminus \Sigma) \ .
\label{stab} \end{equation}

 The structure of the {\em stable} rational cohomology groups of spaces $\PP^k_{d,n}\setminus \Sigma$ with
fixed $k>n$ and growing $d$, that is, of the inverse limits of these groups with respect to the maps (\ref{stab}), follows from Theorem \ref{cor1} by means of a result of \cite{KY}, indicating also a realistic estimate for the instant of the stabilization.
The answers are different for the stabilizations over even and
odd $d$, and are as follows. 
The restrictions of the systems of the class $\PP^k_{d,n}$ 
to the unit sphere in $\R^n$ define injective maps
\begin{equation} \label{embev}\PP^k_{d,n} \setminus \Sigma \to
\mbox{EM}_i(S^{n-1} , (\R^k\setminus 0)) \sim \mbox{EM}_i(S^{n-1}, S^{k-1}),
\end{equation} 
$i \equiv d(\mbox{mod } 2)$;
here $\sim$ means homotopy equivalence.

\begin{proposition}[see \cite{KY}] \label{trivi}
If $k>n$ then the maps $($\ref{embev}$)$ induce isomorphisms of cohomology groups in all dimensions strictly below $(k-n)(d+1)$.
\end{proposition}

\noindent
{\bf Remark.} It follows easily from Weierstrass approximation theorem that
for any natural $k$ and $n$, the inverse limit of groups $H^*(\PP^k_{d,n} \setminus \Sigma)$ over the maps $($\ref{stab}$)$ with even or odd $d$ is naturally isomorphic to $H^*(\mbox{\rm EM}_i(S^{n-1}, S^{k-1}))$, $i \equiv d(\mbox{mod } 2)$, see e.g. Proposition 1 in \cite{dan18}. The previous proposition gives a realistic estimate for the instant of stabilization in the case $k>n$.
\medskip

\begin{example} \rm The space $\PP_{1,3}^k \setminus \Sigma$ is homotopy equivalent to the Stiefel manifold $V_3(\R^k)$. For odd $k>3$ the Poincar\'e polynomial of its rational cohomology group is equal to $(1+t^{k-3})(1+t^{2k-3})$. All elements of this group are stable, i.e. they are induced from certain elements of the group $H^*(\mbox{EM}_1(S^2,S^{k-1}),\Q)$ under the map (\ref{embev}) with $n=3$. Indeed, the number 3 in the first bracket should be considered as $n\equiv m+1$, and in the second one just as 3, see the row \{$m$ even, $M$ even\} of Table \ref{oddd} for $m=2$, $M=k-1$. In the case of even $k>3$ the Poincar\'e polynomial of $H^*(V_3(\R^k))$ is equal to $(1+t^{k-1})(1+t^{2k-5})$; only the generator corresponding to $t^{k-1}$ is stable, see the row \{$m$ even, $M$ odd\} of Table \ref{oddd}. 
 The groups $H^*(\PP_{2,3}^k \setminus \Sigma,\Q)$ are found in \cite{i16};
all of them contain a stable subgroup isomorphic to $\Q$ in dimension $k-1$, see the rows \{$m$ even\} of Table \ref{evv}.
\end{example}

\subsection{Notation} The calculations summarized in Theorems \ref{cor1} and \ref{m4} will be performed in the following terms.

Given a topological space $X$ and local system $L$ of groups on it,
{\bf $\bar H_*(X, L)$} is the (Borel--Moore) homology group of the complex of locally finite singular chains of $X$ with coefficients in $L$.

For any natural $N$, {\bf $I(X,N)$} and {\bf $B(X,N)$} denote respectively the ordered and unordered $N$-point configuration spaces of $X$. In particular, $I(X,N)$ is an open subset in $X^N$, and $B(X,N)=I(X,N)/S(N)$. 
 
{\bf $\pm \Z$} is the ``sign'' local system of groups on $B(X,N)$ which is locally isomorphic to $\Z$, but the loops in $B(X,N)$ act on its fibers as multiplication by $\pm 1$ depending on the parity of the permutations of $N$ points defined by these loops. Also, 
we denote by $\pm \Q$ the local system $\pm \Z \otimes \Q$. 

 If $X$ is a connected manifold, then by $\mbox{Or}$ we denote the local system on $B(X,N)$ which is locally equivalent to $\Z$,
 and moves any element of the fiber to its opposite after the transportation over a loop in $B(X,N)$ if and only if the union of traces of all $N$ points during the movement along this loop defines an element of $H_1(X,\Z_2)$ destroying the orientation of $X$. 

 \begin{lemma} If $X$ is a connected $m$-dimensional manifold, then 
the orientation sheaf of the space $B(X,N)$ is equal to \rm $(\pm \Z)^{\otimes m} \otimes \mbox{Or}.$ \hfill $\Box$
\label{orient}
\end{lemma}

Consider the unique non-trivial local system of groups with fiber $\Z$ on $\RP^m$ : the non-zero element of the group $\pi_1(\RP^m)$ acts on its fibers as multiplication by $-1$. If $m$ is odd then the homology group of $\RP^m$ with coefficients in this system is finite in all dimensions. Denote by $\Theta$ the local system on the space $(\RP^m)^{N}$ (i.e. on the $N$th Cartesian power of $\RP^m$) which is the tensor product of $N$ systems lifted by the standard projections from these systems on the factors. 

Denote by $\tilde \Theta$ the local system on $B(\RP^m,N)$ with fiber $\Z$,
such that a loop in $B(\RP^m,N)$ acts as multiplication by $-1$ on the fiber if and only if
the union of traces of all $N$ points of the configuration during the movement along this loop defines the non-zero element of the group $H_1(\RP^m,\Z_2)$. 

\begin{remark} \label{protriv}
$\tilde \Theta =\mbox{\rm Or}$ if $m$ is even. 
The orientation sheaf of $B(\RP^m,N)$ is equal to $\tilde \Theta$ if $m$ is even, and to $\pm \Z$ if m is odd.
The restriction of the local system $\Theta$ to the subspace $I(\RP^m,N) \subset (\RP^m)^N$ is isomorphic to the system lifted from $\tilde \Theta$ by the standard covering $I(\RP^m,N) \to B(\RP^m,N)$. 
\end{remark}

Denote by $L_r^m$ the lens space $S^m /(\Z/r\Z)$ ($m$ odd). 

$\RP_\star^m$ (respectively, $L^m_{r,\star}$) is the notation for the space $\RP^m$ (respectively, $L^m_{r})$ with one point removed.

\subsection{The main spectral sequence}
\label{agen}

In this subsection and the next one we describe a spectral sequence calculating the cohomology groups of some spaces of equivariant maps, including the ones considered in Theorems \ref{cor1} and \ref{m4}. Its consequences needed for the proof of these two theorems are collected in Corollary \ref{mainss} on page \pageref{mainss}. In \S \ref{which} we discuss which spaces of equivariant maps are homotopy equivalent to ones described here.

Let $X \subset \R^a$ be a compact semialgebraicl set, and $G$ be a compact Lie group acting freely and algebraically on a neighbourhood of $X$ in $\R^a$. Let $\rho:G \to \mbox{O}(\R^W)$ be a representation of $G$, and $\Lambda$ a closed semi-algebraic subset in the unit sphere $S^{W-1} \subset \R^W$, invariant under the corresponding action of $G$ on $S^{W-1}$. Denote by $C\Lambda$ the union of the origin in $\R^W$ and of all rays that start from the origin and intersect the sphere $S^{W-1}$ at the points of $\Lambda$. 

Denote by $\mbox{EM}_G(X,\R^W)$ the space of all continuous maps $X \to \R^W$ equivariant under our actions of $G$ (i.e., the maps $f$ such that $f(g(x))= \rho(g)(f(x))$ for any $x \in X$, $g \in G$). 
Denote by $\mbox{EM}_G(X,(\R^W \setminus C\Lambda))$ its subspace consisting of maps whose images do not meet $C\Lambda$.
If the sum of the dimensions of $X/G$ and $C\Lambda$ is at most $W-2$, then the cohomology group $\tilde H^*(\mbox{EM}_G(X,(\R^W \setminus C\Lambda))$ of this space (reduced modulo a point) can be calculated by a spectral sequence, all whose groups $E^{p,q}_1$ are finitely generated and can be explicitly expressed in the terms of our data. All these non-trivial groups $E^{p,q}_1$ lie in the wedge \begin{equation}
\label{wedge}
\{(p,q)| p < 0, q+p(W-\dim X + \dim G - \dim C\Lambda ) \geq 0\} .
\end{equation} 
In particular, for any $c$ there are only finitely many such non-trivial groups in the line $\{p+q=c\}$. Let us describe these groups.

For any natural $N$ consider the space $\tilde I(X,N)$ of ordered collections of $N$ points of $X$ such that the $G$-orbits of points are pairwise distinct, and consider the trivial fiber bundle over this space with fiber $(C\Lambda)^N$. It is convenient to consider the fiber of this bundle over a sequence $(x_1, \dots, x_N) \in X^N$ as the space of the ways to associate a point $\lambda_j \in C\Lambda$ with any point $x_j$ of this sequence. Denote by ${\mathfrak G}_N $ the semidirect product of the group $G^N$ and the permutation group $S(N)$ with the multiplication $$\left((g_1, \dots, g_N);P\right)\left((h_1, \dots, h_N);Q\right) = \left((g_1h_{P^{-1}(1)}, \dots, g_Nh_{P^{-1}(N)});PQ\right)$$
for any $g_j, h_j \in G$ and $P, Q \in S(N)$. This group acts on the base and the fibers of our fiber bundle: its element 
\begin{equation}
\label{semi}
((g_1, \dots, g_N); P)
\end{equation} with $g_j \in G$, $P \in S(N)$, takes the points $(x_1, \dots, x_N)$ and $(\lambda_1, \dots, \lambda_N)$ to 
\begin{equation} \label{semact}
\left(g_1(x_{P^{-1}(1)}), \dots, g_N(x_{P^{-1}(N)})\right) \qquad \mbox{and} \qquad \left(\rho(g_1)(\lambda_{P^{-1}(1)}), \dots, \rho(g_N)(\lambda_{P^{-1}(N)})\right)
\end{equation} respectively. Therefore this group acts also on the space of this fiber bundle. Denote by $\Pi_{N,G}(X)$ the quotient space of the space of our fiber bundle by this action; it can be considered as the space of a fiber bundle with the same fiber $(C\Lambda)^N$ over the configuration space $B(X/G,N)$. 

Also, define the {\em sign} of the element 
(\ref{semi}) (and of entire its connected component in ${\mathfrak G}_N$) as the product of the determinants of all $N$ operators $\rho(g_j) \in \mbox{O}(\R^W)$ and the $(W+1)$th power of the sign of the permutation $P$. Consider the homomorphism $\pi_1(B/G,N) \to \pi_0({\mathfrak G}_N)$ from the exact sequence of the principal ${\mathfrak G}_N$-bundle $\tilde I(X,N) \to B(X/G,N)$. Composing it with the sign homomorphism $\pi_0({\mathfrak G}_N) \to \{\pm 1\}$ we obtain a homomorphism $\pi_1(B/G,N) \to \{\pm 1\}$. We denote the resulting local system with fiber $\Z$ by ${\mathcal J}$, and we use the same notation ${\mathcal J}$ for the pullback of this local system to $\Pi_{N,G}(X)$.

\begin{theorem}
\label{absv}
If the sum of dimensions of spaces $X/G$ and $C\Lambda$ does not exceed $W-2$, then there is a spectral sequence $E^{p,q}_r$ converging to the reduced homology group $\tilde H^*( \mbox{\rm EM}_G(X,(\R^W \setminus C\Lambda)))$; all its groups $E^{p,q}_1$ with $p\geq 0$ are trivial, and 
\begin{equation}
\label{absvf}
E^{p,q}_1 \simeq \bar H_{-pW-q}(\Pi_{-p,G}(X), {\mathcal J}) 
\end{equation}
for $p<0$.
\end{theorem} 

In particular, for any $p<0$ the groups $E^{p,q}_1$ with this $p$ are trivial if $q$ does not belong to the segment $[-p(W-\dim X +\dim G - \dim C\Lambda), -pW]$. Indeed, only for such $q$ the lower index $-pW-q$ in (\ref{absvf}) belongs to the segment $[0, \dim \Pi_{-p,G}(X)]$.

\subsection{The relative version}

Let in conditions of Theorem \ref{absv} additionally $A $ be a $G$-invariant subcomplex in $X$, and $\varphi: A \to \R^W \setminus C\Lambda$ be a $G$-equivariant map. Denote by $\mbox{EM}_G(X,A;\R^W)$ the space of all $G$-equivariant maps $X \to \R^W$ coinciding with $\varphi$ on $A$. Denote by $\Pi_{N,G}(X,A)$ the part of $\Pi_{N,G}(X)$ consisting of all fibers 
of our fiber bundle $\Pi_{N,G}(X) \to B(X/G,N)$ over the points of the subspace $B((X\setminus A)/G,N) \subset B(X/G,N)$. 

\begin{theorem}
\label{relav}
If the sum of dimensions of spaces $X/G$ and $C\Lambda$ does not exceed $W-2$, then there is a spectral sequence $\hat E^{p,q}_r$ converging to the reduced homology group $\tilde H^*( \mbox{\rm EM}_G(X,A;(\R^W \setminus C\Lambda))$; all its groups $\hat E^{p,q}_1$ with $p\geq 0$ are trivial, and 
\begin{equation}
\label{relaf}
\hat E^{p,q}_1 \simeq \bar H_{-pW-q}(\Pi_{-p,G}(X,A), {\mathcal J}) 
\end{equation}
for $p<0$.
\end{theorem} 

\begin{remark}\rm
In the case of the trivial group $G$, the spectral sequences (\ref{absvf}) and (\ref{relaf}) were
 defined in \S III.6 of \cite{Book}. The sequence (\ref{absvf}) is in this case conjecturally equipotent to the one by Anderson \cite{anderson}: they are formally non-isomorphic since the first invariant term of the Anderson's spectral sequence is $E_2$ and that of ours one is $E_1$, but probably they handle effectively nearly the same class of examples. 
The construction of spectral sequences (\ref{absvf}), (\ref{relaf}) is bery similar to the one given in \cite{Book}; we review it with appropriate modifications in \S \ref{sseven}. 
\end{remark}

\begin{corollary}
\label{mainss}
For any $m<M$, there are three spectral sequences $\{E^{p,q}_r\}$ calculating the integer cohomology groups of the spaces $\mbox{Map}(S^m,S^M)$, $\mbox{EM}_0(S^m,S^M)$, and $\mbox{EM}_1(S^m,S^M)$, and three more spectral sequences $\{\hat E^{p,q}_r\}$ for the cohomology of spaces $\Omega^m S^M$, $\mbox{EM}^*_0(S^m,S^M)$, and $\mbox{EM}^*_1(S^m,S^M)$; 
the supports of all six spectral sequences lie in the wedge $\{(p,q): p \leq 0$, $q+p(M-m+1) \geq 0\}$, the unique non-trivial groups $E^{p,q}_1$ and $\hat E^{p,q}_1$ with $p=0$ are $E^{0,0}_1 \simeq \hat E^{0,0}_1 \simeq \Z$; the group $E_1^{p,q}$ with $p <0$ is equal to 
\begin{eqnarray}
\bar H_{-p(M+1)-q}\left(B(S^m, -p), (\pm \Z)^{\otimes M}\right) \qquad \qquad & \mbox{for } &\mbox{Map}(S^m,S^M) \ , 
\label{msgen} \\
\label{ssss}
\bar H_{-p(M+1)-q}\left(B(\RP^{m},-p),(\pm \Z)^{\otimes M}\right) \qquad \qquad & \mbox{for } & \mbox{EM}_0(S^m,S^M), \\
\label{thomm2}
\bar H_{-p(M+1)-q}\left(B(\RP^{m},-p),(\pm \Z)^{\otimes M} \otimes \tilde \Theta^{\otimes (M+1)}\right) & \mbox{for } & \mbox{EM}_1(S^m,S^M) \ ;
\end{eqnarray}
the analogous groups $\hat E^{p,q}_1$ with $p<0$ are equal to 
\begin{eqnarray}
\label{msgenpt}
\bar H_{-p(M+1)-q}\left(B(\R^m, -p), (\pm \Z)^{\otimes M}\right) \qquad \qquad & \mbox{for } & \Omega^m S^M, \\
\label{sssspt}
\bar H_{-p(M+1)-q}\left(B(\RP_\star^{m},-p),(\pm \Z)^{\otimes M}\right) \qquad \qquad & \mbox{for } & \mbox{EM}^*_0(S^m,S^M) \ , \\
\label{thomm2pt}
\bar H_{-p(M+1)-q}\left(B(\RP_\star^{m},-p),(\pm \Z)^{\otimes M} \otimes \tilde \Theta^{\otimes (M+1)}\right) & \mbox{for } & \mbox{EM}^*_1(S^m,S^M) .
\end{eqnarray}

For any odd numbers $m<M$, there is a spectral sequence $\{E^{p,q}_r\}$ $($respectively, $\{\hat E^{p,q}_r\})$ calculating the integer cohomology group of the space of maps considered in the first $($respectively, the second$)$ statement of Theorem \ref{m4}; their groups $E_1^{p,q}$ and $\hat E_1^{p,q}$ for any $p <0$ are equal respectively to 
\begin{equation}
\label{lenspm}
\bar H_{-p(M+1)-q}\left(B(L_r^{m},-p),\pm \Z\right) \qquad \mbox{and} 
\end{equation}
\begin{equation}
\label{lenspmpt}
\bar H_{-p(M+1)-q}\left(B(L_{r,\star}^{m}\ ,-p),\pm \Z\right) \ .
\end{equation}
\end{corollary}

\noindent
{\it Proof.} The formulas (\ref{msgen}--\ref{thomm2}) and (\ref{lenspm}) follow immediately from Theorem \ref{absv}, and formulas (\ref{msgenpt}--\ref{thomm2pt}) and (\ref{lenspmpt}) from Theorem \ref{relav}. In all these cases $W=M+1$ and $\Lambda=\emptyset$ (so that $C\Lambda$ is the origin in $\R^W$). In (\ref{ssss}) and (\ref{sssspt}) $X=\RP^m$
(since the even maps $S^m \to S^M$ are in the obvious one-to-one correspondence with usual maps $\RP^m \to S^M$), and group $G$ is trivial. Analogously, the maps $S^m \to S^M$ commuting with the actions of the group $\Z_r$ according to the rule (\ref{chara}) are in one-to-one correspondence with the maps $L_\tau^m \to S^M$ (where $\tau$ is the greatest common divisor of $r$ and $s$) commuting with some free actions of the group $G=\Z_\frac{r}{\tau}$: this action on $L^m_\tau$ is the factorization of the action on $S^m$. Therefore in (\ref{lenspm}) and (\ref{lenspmpt}) we have $X = L^m_\tau$. In all other cases $X=S^m$: the group $G$ equals $\Z_2$ in (\ref{thomm2}) and (\ref{thomm2pt}), and is trivial in (\ref{msgen}) and (\ref{msgenpt}).
In these examples it is convenient for us to consider the non-reduced homology groups, therefore we add the standard cells $E^{0,0} \simeq \Q$ in comparison to the spectral sequences of Theorems \ref{absv} and \ref{relav}.
\hfill $\Box$ 

\begin{remark} \rm 
Spectral sequence (\ref{msgenpt}) stabilizes at the term $E_1$ by the Snaith splitting \cite{snaith}. The calculations given below prove a similar stabilization result of all spectral sequences (\ref{msgen}--\ref{lenspmpt}) at least over the rational coefficients.
Formulas (\ref{msgen}--\ref{ssss}) and (\ref{msgenpt}--\ref{sssspt}) related with the trivial group actions are the special cases of the main formula of \S III.6.2 in \cite{Book}. Formula (\ref{thomm2}) is announced in \cite{dan18}, together wit the ``non-pointed'' part of Theorem \ref{cor1} concerning the third columns of Tables \ref{evv} and \ref{oddd}. 
\end{remark}

\subsection{Twisted homology of configuration spaces}
In order to deduce Theorems 1 and 2 from these spectral sequences, we need to calculate the homology of several configuration spaces, in particular the following ones.

\begin{theorem}
\label{lem29}
For any $r \geq 1$ and any odd $m$,

$($a$)$ $H^*(B(L_r^{m},N),\Q) \simeq H^*(S^{m},\Q)$ for all $N \geq 1$, 

$($b$)$ the group $H^j(B(L_r^{m},N),\pm \Q)$ is
isomorphic to $\Q$ if $N$ is odd and $j$ is equal to either $(m-1)\frac{N-1}{2}$ or $(m-1)\frac{N-1}{2}+m,$
and is trivial for all other combinations of $N$ and $j$.
\end{theorem}

\begin{theorem} \label{main2}
Let $m$ be odd. Then 

A$)$ The group $H_j(B(\RP_\star^m,N), \tilde \Theta \otimes \pm \Q)$ with arbitrary $N \geq 1$ is equal to $\Q$ if $j= \left]\frac{N}{2}\right[ \times (m-1)$ $($where $]X[$ denotes the smallest integer not smaller than $X)$ and is trivial for all other $j$. 

B$)$ The $($Poincar\'e dual to one another$)$ groups $H_*(B(\RP^m,N), \tilde \Theta \otimes \pm \Q)$ and $\bar H_*(B(\RP^m,N), \tilde \Theta \otimes \Q)$ are trivial if $N$ is odd, and have the Poincar\'e polynomials equal respectively to $t^{\frac{N}{2}(m-1)}(1+t^m)$ and $ t^{\frac{N}{2}(m+1)} (1+t^{-m})$ if $N$ is even.
\end{theorem}

\begin{theorem} \label{main1}
If $m$ and $N$ are odd, then the group $H_*(I(\RP^m,N), \Theta \otimes \Q)$ is trivial.
\end{theorem}

The next theorem is not needed for our calculations, nevertheless I prove it below for the completeness.

\begin{theorem} \label{lem1} If $N$ is even and $m$ is odd, then
the Poincar\'e polynomial of the group $H_*(I(\RP_\star^m,N-1), \Theta \otimes \Q)$ is equal to
\begin{equation} \label{fibfor}
(N-1)!! \ t^{(m-1)N/2} 
\prod_{r=1}^{N/2-1} \left(1+(2r-1)t^{m-1}\right) .
\end{equation}
\end{theorem}

\section{Proof of Theorem \ref{lem29}}

\begin{lemma}
\label{lem11}
For any $N \geq 1$, an arbitrary embedding $\R^m \to L_{r,\star}^m$ induces the isomorphisms 
\begin{eqnarray} 
 \label{a}
H^*(I(L_{r,\star}^m,N),\Q) & \simeq & H^*(I(\R^{m},N),\Q); \\
\label{b}
H^*(B(L_{r,\star}^m, N), \Q) & \simeq & H^*(B(\R^{m}, N), \Q); \\
\label{c}
H^*(B(L_{r,\star}^m, N), \pm \Q) & \simeq & H^*(B(\R^{m}, N), \pm \Q).
\end{eqnarray}
\end{lemma}

\noindent
{\it Proof.} The isomorphism (\ref{a}) follows by induction over the standard fiber bundles $I(X,l+1) \to I(X,l)$ defined by forgetting the last points of configurations, with $X$ equal to either $L^m_{r,\star}$ or $\R^m$. Indeed, 
we have natural isomorphisms $ H^*(\R^m \setminus \{l \ \mbox{ points}\},\Q) \simeq 
H^*( L_{r,\star}^m \setminus \{l \ \mbox{ points}\},\Q)$ for their fibers. These bundles are homologically trivial under the action of $\pi_1(I(X,l))$ since the cohomology groups of fibers are generated by the linking number classes with the removed ordered points. Therefore we also have natural isomorphisms of the spectral sequences of these fiber bundles. 

These spectral sequences for $I(\R^m,N)$ stabilize at the terms $E_2$, so that the Poincar\'e polynomials of the groups (\ref{a}) are equal to \begin{equation}
\label{fks}
\prod_{a=1}^{N-1}(1+ a t^{m-1}) \ ,
\end{equation}
 cf. \cite{A69}, \cite{Cohen76}.
 
By Proposition 2 from \cite{Serre}, any group considered in (\ref{b}) or (\ref{c}) is equal to the quotient of the group considered in the same side of (\ref{a}) by the subgroup generated by all elements of the form $g(\alpha)-\alpha$, where $g$ is an element of the permutation group $S(N)$ that acts as follows on the cohomology of the spaces $I(L_{r,\star}^m,N)$ and $I(\R^m,N)$. In the case of coefficients in $\Q$ considered in (\ref{b}) this action of $S(N)$ is the obvious one induced by the action of the permutations on the fibers of a principal $S(N)$-covering; in the case of $\pm \Q$-coefficients any element of $S(N)$ additionally multiplies the classes by $\pm 1$ depending on the parity of the corresponding permutation. 
These $S(N)$-actions also commute with the cohomology maps induced by the embedding $I(\R^m,N) \to I(L_{r,\star}^m,N)$, hence the results of this factorization for $L_{r,\star}^m$ are the same as for $\R^m$. \hfill $\Box$ \medskip

\begin{proposition}[see \cite{Cohen76}]
\label{chn76}
For any odd $m$, both groups $ H^*(B(\R^{m}, N), \Q)$ and \\ $H^*(B(\R^{m}, N), \pm \Q)$ are one-dimensional over $\Q$. Namely, 
$ H^0(B(\R^{m}, N), \Q)\simeq \Q$ and $H^{(m-1)[N/2]}(B(\R^{m}, N), \pm \Q) \simeq \Q$.
\end{proposition}

The first assertion of this proposition holds by the Euler characteristic reason: by (\ref{fks}) the cohomology group $H^*(I(\R^{m},N), \Q)$ of the $N!$-fold covering space is exactly $N!$-dimensional over $\Q$ and is placed in even dimensions only. 

The groups (\ref{b}) are obviously isomorphic to $\Q$ in dimension $0$, hence they are trivial in all other dimensions. 

Let us describe a cycle generating the homology group dual to (\ref{c}).
Let $N$ be even. Fix $N/2$ different points $a_s \in \R^m$, $s=1, \dots, N/2$, and a small number $\varepsilon>0$. Consider the submanifold $V(N) \subset B(\R^m,N)$ consisting of all configurations of $N$ points split into $N/2$ pairs, the points of the $s$th pair being some opposite points of the $\varepsilon$-sphere around $a_s$ in $\R^m$. If $N$ is odd then a similar $((m-1)[N/2])$-dimensional manifold $V(N)$ is obtained from the manifold $V(N-1)$ we have just defined by adding one constant point $X_N$ distant from all points $a_s$ to all $(N-1)$-configurations forming this manifold $V(N-1)$.

\begin{lemma} \label{lemW} If $m$ is odd then the manifold $V(N)$ is $\pm \Z$-orientable for any $N$, and
the group $H_*(B(\R^m,N),\pm \Q)$ is generated by its fundamental $\pm \Q$-cycle.
\end{lemma}

{\it Proof} (cf. \cite{fuchs}). The assertion on the $\pm \Q$-orientability is immediate. Let us fix a direction in $\R^m$. If $N$ is even, denote by $\Xi$ the subset in $ B(\R^m,N)$ consisting of all configurations, whose $N$ points can be split into pairs, such that the line through any pair has the chosen direction. For odd $N$, let $\Xi$ consist of all such $(N-1)$-configurations augmented arbitrarily with some $N$th points. $\Xi$ is an $(mN - (m-1)[N/2])$-dimensional semi-algebraic variety with a standard orientation in its regular points, and its fundamental cycle defines an element of the group 
$\bar H_{mN - (m-1)[N/2]}(B(\R^m,N), \Q)$. In particular, its intersection index with the cycle $V(N) \subset B(\R^m,N)$ is well defined. If the points $a_s$ (and $X_N$ if $N$ is odd) defining $W$ are in general position with respect to the chosen direction, and $\varepsilon$ is small, then these two cycles have exactly one transversal intersection point, therefore the homology classes of both of them are not equal to 0. \hfill $\Box$ \medskip

Consider now the exact sequence of the pair $(B(L^m_{r},N), B(L^m_{r,\star},N))$ with coefficients in the system $\pm \Q$. 
Its term $H_i(B(L^m_{r},N), B(L^m_{r,\star},N); \pm \Q)$ is isomorphic to 
$H_{i-m}(B(L^m_{r,\star},N-1), \pm \Q) $ by the Thom isomorphism of the (trivial) normal bundle of the space $B(L^m_{r},N) \setminus B(L^m_{r,\star},N) \simeq B(L^m_{r,\star},N-1) $ in $B(L^m_{r},N).$ By formula (\ref{c}) and Proposition \ref{chn76} this exact sequence contains only two non-zero (and equal to $\Q$) terms not of the form $H_i(B(L^m_{r},N), \pm \Q)$. 
 In the case of odd $N$ these two terms are distant from one another in this sequence and imply Theorem \ref{lem29} for such $N$. If $N$ is even then these two terms form the fragment 
\begin{equation}H_{i+1}(B(L^m_{r},N), B(L^m_{r,\star},N); \pm \Q) \to H_{i}(B(L^m_{r,\star},N), \pm \Q)\label{ii}
\end{equation} 
 with $i=(m-1)N/2.$ 
The image of the basic cycle of the left-hand group under this map can be realized by the $(m-1)\frac{N}{2}$-dimensional manifold in $B(L^m_{r,*},N)$ swept out by all $N$--configurations, in which some $N-1$ points form any configurations from the manifold $V(N-1)$ (in particular one of these points is fixed), and one point more belongs to a small $(m-1)$-dimensional sphere around the point $L^m_{r}\setminus L^m_{r,\star}$. This cycle is homologous to a similar one, in the definition of which the last sphere is replaced by a huge sphere in $\R^m \subset L^m_{r,\star} $ encircling all other points of all configurations of our cycle. The intersection index of the obtained cycle with the cycle $\Xi$ from the proof of Lemma \ref{lemW} is equal to 2 or $-2$ depending on the choice of orientations, therefore the map (\ref{ii}) is an isomorphism, and all groups $H_i(B(L^m_{r},N), \pm \Q)$ are trivial.\hfill $\Box$

\begin{remark} \rm Part A) of Theorem \ref{lem29} can also be deduced by applying the methods of the paper \cite{BCL} which contains a general algorithm for calculating the Betti numbers of configuration spaces of odd-dimensional manifolds. 
\end{remark}

\section{Proof of Theorems \ref{main2}, \ref{main1}, and \ref{lem1}}

\subsection{Proof and realization of Theorem \ref{main2}(A)}

\begin{lemma}
\label{lempodd}
If $m$ is odd, then the group $H_*(I(\RP^m_\star,N), \Theta \otimes \Q)$ with any $N$ is $N!$-dimensional over $\Q$, and the group $H_*(B(\RP^m_\star,N), \tilde \Theta \otimes \pm \Q)$ is one-dimensional. All groups $H_j(I(\RP^m_\star,N), \Theta \otimes \Q)$ and $H_j(B(\RP^m_\star,N), \tilde \Theta \otimes \pm \Q)$ with odd $j$ are trivial.
\end{lemma}

\noindent
{\it Proof.} There is a principal $(\Z_2)^N$-covering over the space $I(\RP^m_\star,N)$: it is the space of all sequences of $N$ different points in $S^m$ such that none of its elements coincides with either of the two preimages of the distinguished point of $\RP^m$ under the standard covering map $S^m \to \RP^m$, or with a point opposite in $S^m$ to an other element of the sequence. 

The space $X_N$ of this covering has the structure of a tower of fiber bundles $X_N \to X_{N-1} \to \dots \to X_2 \to X_1$ where $X_1$ is the sphere $S^m$ with $2$ points removed, and the fiber of any map $X_{k} \to X_{k-1}$ is $S^m$ with $2k$ points removed. It follows easily from the spectral sequences of all these bundles that the Poincar\'e polynomial of the rational homology group of the total space of this tower is equal to $\prod_{j=1}^{N} (1+(2j-1)t^{m-1})$. In particular, all its homology groups lie in even dimensions only, and its Euler characteristic is equal to $2^N N!$. 

The group $(\Z_2)^N$ of this covering acts on the complex of $\Q$-singular chains of the space of the covering: an element $g \in (\Z_2)^N$ transforms any singular simplex geometrically via the usual action of this group by central symmetries on the factors $S^m$ of the space $(S^m)^N$ containing our covering space, and additionally multiplies it by $-1$ if the number of ones in $g$ is odd. The group $H_*(I(\RP^m_\star,N), \Theta \otimes \Q)$ can be considered as the homology group of the complex of coinvariants for this action, hence this group also has non-trivial elements in even dimensions only. In particular, the total dimension of this group is equal to its Euler characteristic, which is equal to the Euler characteristic $2^N N!$ of the space of the covering divided by the degree $2^N$ of this covering. 

In a similar way, $H_*(B(\RP^m_\star,N), \tilde \Theta \otimes \pm \Q)$ is the homology group of the complex of coinvariants for an action of a group of order $2^N N!$ on rational singular chains, hence it also has no non-trivial odd-dimensional elements, and its total dimension is equal to its Euler characteristic and hence to 1. \hfill $\Box$ \medskip 

For even $N$, let $V(N) \subset B(\R^m,N) \subset B(\RP^m_\star,N)$ be the manifold from the proof of Lemma \ref{lemW} (we assume that all points $a_s$ participating in its construction are distant from the removed point $\RP^m \setminus \RP^m_\star$). For odd $N$, consider the $\left(\, \left]\frac{N}{2}\right[\times (m-1)\right)$-dimensional manifold $U(N)$ in $B(\RP^m_\star,N)$ the elements of which are unions of an $N-1$-configuration from $V(N-1) \subset B(\RP^m_\star,N-1)$ and a one point subset of the $\varepsilon$-sphere about the point $\RP^m \setminus \RP^m_\star$.

\begin{lemma}
\label{realiz}
The basic cycle of the manifold $V(N)$ $($in the case of even $N)$ or $U(N)$ $($if $N$ is odd$)$ defines a non-trivial element of the group $H_*(B(\RP^m_\star,N), \tilde \Theta \otimes \pm \Q)$.
\end{lemma}

\noindent
{\it Proof.} Let us describe a Poincar\'e dual element of the group $\bar H_*(B(\RP^m_\star,N), \tilde \Theta \otimes \Q)$. 

For even $N$, consider the $\frac{N}{2}(m+1)$-dimensional subvariety in $B(\RP_\star^m,N)$ consisting of all $N$-configurations which can be split into pairs so that the elements of each pair belong to the same fiber of the Hopf fibration $\RP^m \to \CP^{\frac{m-1}{2}}$. Its regular locus is $\tilde \Theta$-orientable, and it defines an element of the group $\bar H_{\frac{N}{2}(m+1)}(B(\RP_\star^m,N),\tilde \Theta)$.

 If $N$ is odd, then we consider the space of all $N$-configurations, which are the unions of an $N-1$-element subset as in the previous paragraph, and a one point subset of the fiber of the Hopf fibration through the point $\RP^m \setminus \RP^m_\star$. 

The intersection indices of cycles $V(N)$ or $U(N)$ (depending on the parity of $N$) with these cycles are equal respectively to $\pm 1$ and $\pm 2$, in particular the classes of all these cycles in the corresponding homology groups are not equal to 0. \hfill $\Box$ \medskip

Theorem \ref{main2}(A) is thus also proved. \hfill $\Box$

\subsection{Simplicial resolution of diagonal varieties \rm (see \cite{noneven})}
\label{sire}

Denote by $\nabla$ the set $(\RP^m)^N \setminus I(\RP^m,N)$. It is a stratified variety, whose strata correspond to all partitions of the set $\{1, \dots, N\}$ into $\leq N-1$ non-empty parts. Namely, with any such partition $A$ into some sets $A_1, \dots, A_k$ the closed stratum $L(A) \subset \nabla$ is associated, which consists of all collections $(x_1, \dots, x_N) \in (\RP^m)^N$ such that $x_\alpha =x_\beta$ whenever $\alpha$ and $\beta$ belong to the same part of this partition. Obviously, $L(A) \simeq (\RP^m)^k$.

Such partitions form a partially ordered set: any partition dominates its subpartitions. Denote by $\Delta(N)$ the order complex of this poset.

Let $m$ be odd, then by the Poincar\'e--Lefschetz duality in the $(\Theta \otimes \Q)$-acyclic space $(\RP^m)^N$ we have 
\begin{equation} \label{lefsch}
H^i(I(\RP^m,N), \Theta \otimes \Q) \simeq H_{mN-i-1}(\nabla, \Theta \otimes \Q)
\end{equation}
for any $i$,
so we will consider the right-hand groups in (\ref{lefsch}). We use the standard simplicial resolution ${\mathfrak S}$ of $\nabla$, which is a subset in $\Delta(N) \times \nabla$ defined as follows. For any vertex of the complex $\Delta(N)$ (i.e., some non-complete partition $A$) denote by $\Delta_A \subset \Delta(N)$ its subordinate subcomplex, i.e. the order complex of subpartitions of $A$ (including $A$ itself). Let $\partial \Delta_A$ be the link of $\Delta_A$, that is, the union of its simplices not containing the maximal vertex $\{A\}$. The space ${\mathfrak S} \subset \Delta(N) \times \nabla$ is the union of the sets $\Delta_A \times L(A)$ over all non-complete partitions $A$ of $\{1, \dots, N\}$. Let $\Theta^!$ be the local system on ${\mathfrak S}$ lifted from $\Theta$ along the canonical projection ${\mathfrak S} \to \nabla$. We have the standard isomorphism
\begin{equation}
\label{simresi}
H_* ({\mathfrak S}, \Theta^!) \simeq H_*(\nabla, \Theta). 
\end{equation}

The space ${\mathfrak S}$ has a natural filtration ${\mathfrak S}_0 \subset \dots \subset {\mathfrak S}_{N-2} = {\mathfrak S}$: its term ${\mathfrak S}_p$ is the union of spaces $\Delta_A \times L(A)$ over all partitions of $\{1, \dots,N\}$ into $\geq N-1-p$ sets. Consider the spectral sequence calculating the group $H_*({\mathfrak S}, \Theta^! \otimes \Q)$ and induced by this filtration.
Any set ${\mathfrak S}_p \setminus {\mathfrak S}_{p-1}$ splits into the blocks corresponding to all partitions $A$ into exactly $N-1-p$ parts and equal to 
\begin{equation}
\label{block}
(\Delta_A \setminus \partial \Delta_A) \times L(A).
\end{equation}
 Correspondingly, the term $E^1_{p,q}$ of this spectral sequence splits into the sum of groups $
\bar H_{p+q}((\Delta_A \setminus \partial \Delta_A) \times L(A), \Theta^!)$
 over such partitions $A$.

\begin{lemma}
The group $\bar H_*((\Delta_A \setminus \partial \Delta_A) \times L(A),\Theta^! \otimes \Q)$ 
 is trivial if some parts of $A$ consist of odd numbers of elements. If all parts of $A$ are even, then this group is equal to $\bar H_*(\Delta_A \setminus \partial \Delta_A, \Q) \otimes H_*((\RP^m)^{k(A)}, \Q)$, where $k(A)$ is the number of parts of $A$. 
\end{lemma} 

This lemma follows immediately from the K\"unneth formula, see \cite{noneven}. \hfill $\Box$ \medskip

Theorem \ref{main1} follows immediately from this lemma and the isomorphisms (\ref{lefsch}) and (\ref{simresi}). \hfill $\Box$

\subsection{Proof of Theorem \ref{lem1}.}

The projection of $(\RP^m)^N$ onto its first factor defines a fiber bundle $\nabla \to \RP^m$. Let $\tilde \nabla$ be its fiber over a fixed point $x_1$. The space $I(\RP_\star^m,N-1)$ is equal to its complement $(\RP^m)^{N-1} \setminus \tilde \nabla$ in the product of remaining $N-1$ factors. 
The restriction of our simplicial resolution to the fiber $\tilde \nabla$ defines a simplicial resolution $\tilde {\mathfrak S}$ of $\tilde \nabla$. Again, we have the Lefschetz isomorphism 
\begin{equation}
\label{lef2}
H^i(I(\RP_\star^m,N-1), \Theta \otimes \Q) \simeq H_{m(N-1)-i-1}(\tilde \nabla, \Theta \otimes \Q) \simeq H_{m(N-1)-i-1}(\tilde {\mathfrak S}, \Theta^! \otimes \Q)
\end{equation}
for any $i$; let us calculate the right-hand groups in this equality using the restriction to $\tilde {\mathfrak S}$ of our filtration on ${\mathfrak S}$.
The term $\tilde E^1$ of the corresponding spectral sequence $\tilde E^r_{p,q}$ calculating $H_*(\tilde {\mathfrak S}, \Theta^! \otimes \Q)$ was considered in \cite{noneven}, in particular it was shown there that this term has non-trivial groups $\tilde E^1_{p,q}$ only on $N/2$ horizontal segments corresponding to $q=0, m, 2m, \dots, (N/2-1)m$. Namely, for any $s=0,1, \dots, N/2-1$, only the groups $\tilde E^2_{p,sm}$ with natural $p \in [N/2-1, N-2-s]$ are non-trivial. The differential $d_1$ of the spectral sequence supplies these horizontal segments with the structure of chain complexes.

\begin{lemma}
All these horizontal complexes are acyclic in all their terms except for the top-dimensional ones $($that is, the groups $\tilde E^2_{p,q}$ can be non-trivial only for $(p,q) = (N-2,0), (N-3,m), (N-4, 2m), \dots, (N/2-1,(N/2-1)m)\ )$. 
\end{lemma}

\noindent
{\it Proof.} For $m> N/2-1$ this is Lemma 1 in \cite{noneven}; but by the construction these $N/2$ complexes do not depend on $m$ (provided $m$ is odd). \hfill $\Box$

\begin{corollary}
Our spectral sequence stabilizes in the term $\tilde E^2 \equiv \tilde E^\infty.$ The dimensions of its non-trivial groups $\tilde E^2_{p,sm}$ are equal to $\pm$ the Euler characteristic of corresponding horizontal complexes. \hfill $\Box$
\end{corollary}

This corollary reduces Theorem \ref{lem1} to the following lemma.

\begin{lemma} \label{gen}
The Euler characteristic of the horizontal complex corresponding to $q =s m$, $s \in [0, \dots, N/2-1],$ is equal to $(-1)^s (N-1)!!$ times the coefficient at the term $\tau^{N/2-1-s}$ in the polynomial $(1+\tau)(1+3\tau) \cdots (1+(N-3)\tau)$.
\end{lemma}
 
\noindent
{\it Proof.} By the description of the term $\tilde E^1$ given in \cite{noneven}, the dimension of the group $\tilde E^1_{p, sm}$, $p \in [N/2-1, N-2-s]$, can be calculated as follows. Consider all unordered partitions of the set $\{1, \dots, N\}$ into $N-p-1$ parts of even cardinalities, then count any such partition into parts of cardinalities $N_1, \dots, N_{N-p-1}$ with the weight 
\begin{equation} \label{weight}
(N_1-1)!(N_2-1)! \cdots (N_{N-p-1} -1)! \end{equation} 
and multiply the obtained sum by $\binom{N-p-2}{s}.$

Given such a partition, let us order its parts lexicographically: the first part contains the element $1 \in \{1, \dots, N\}$, the second part contains the smallest element not contained in the first part, etc. Then our coefficient (\ref{weight}) is equal to the number of permutations $(a_1, \dots, a_N)$ of $\{1, \dots, N\}$, starting with 1, then running arbitrarily all other elements of the first part, then choosing the smallest element of the second part, then running arbitrarily the remaining elements of this part, etc. The common factor $\binom{N-p-2}{s}$ can be understood as the number of possible selections of some $s$ parts of our partition, not equal to the first one (or, which is the same, of the smallest elements of these parts). 

So, the dimension of the group $E^1_{p,sm}$ is equal to the number of permutations $(a_1, \dots, a_N)$ of $\{1, \dots, N\}$ having some additional furniture and satisfying some conditions. The furniture consists of the choice of two subsets $R \supset S$ 
of cardinalities $N-p-2$ and $s$ respectively in the set of odd numbers $\{3, 5, \dots, N-1\}$; the conditions claim that $a_1=1$ and the numbers $a_i$ with $i \in R$ should be smaller than all forthcoming elements $a_j$, $j >i$.

Consider now an arbitrary permutation of $\{1, \dots, N\}$ with the first element $a_1=1$ and $s$ other selected elements $a_i$ on some odd places 
\begin{equation} \label{places}
i= N_1+1, N_1+N_2+1, \dots, N_1+\dots + N_{s}+1,
\end{equation}
 and count its occurrences from the above construction for all $p \in [N/2-1, N-2-s]$. 
Any of these $s$ selected elements $a_i$ should be smaller than all the forthcoming ones, otherwise our permutation does not occur at all. Further, let $r$ be the number of non-selected and not equal to 1 elements $a_j$ of this permutation with odd $j$, which are smaller than all the forthcoming elements. If $r=0$ then our permutation occurs only once from the consideration of the group $\tilde E^1_{N-2-s,sm}$. If however $r>0$ then it occurs once in $\tilde E^1_{N-2-s, sm}$, $r$ times in $\tilde E^1_{N-3-s,sm}$, $\binom{r}{2}$ times in $\tilde E^1_{N-4-s,sm},$ etc. These $2^r$ occurrences are counted with the alternating signs in the calculation of the Euler characteristic of our horizontal complex, and annihilate one another for $r>0$. Therefore only the permutations with $r=0$ contribute to this Euler characteristic. The number of such permutations (that is, permutations starting with 1, having on the selected $s$ odd places (\ref{places}) $s$ elements which are smaller than all the forthcoming ones, and {\em not} having elements with this property on other odd places) is obviously equal to 
\begin{equation} \label{combi}
(N-1)!!\frac{(N-3)!!}{(N-N_1-1)(N-(N_1+N_2)-1) \cdots (N-(N_1+\dots + N_s)-1)}.
\end{equation}
But the fraction in (\ref{combi}) is a summand after opening the brackets in the product $(1+\tau)(1+3\tau) \cdots (1+(N-3)\tau)$, containing exactly $s$ factors equal to $1$ (coming from the brackets No. $\frac{N-N_1}{2}, \frac{N-(N_1+ N_2)}{2},$\dots, $\frac{N-(N_1+ \dots +N_s)}{2}$ ). \hfill $\Box$

\smallskip
\noindent
{\bf Problem.} Realize a basis in $H_*(I(\RP_\star^m,N-1), \Theta \otimes \Q)$ by explicitly defined cycles.

\subsection{Proof of Theorem \ref{main2}(B)} 
The permutation group $S(N)$ acts in an obvious way on the space $I(\RP^m,N)$ and on the local system $\Theta$ on it. Consider the resulting action of $S(N)$
on the complex of singular $\Theta$-chains of $I(\RP^m,N)$, moving any simplex geometrically by the previous action and additionally multiplying it by $\pm 1$ depending on the parities of the permutations. The group $H_*(B(\RP^m,N), \tilde \Theta \otimes \pm \Q)$ is the homology group of the complex of coinvariants for this action. Again by Proposition 2 from \cite{Serre}, it is also the quotient group of the group 
$H_*(I(\RP^m,N),\Theta \otimes \Q)$ by the subgroup generated by all elements of the form $g(\alpha)-\alpha$ where 
$\alpha \in H_*(I(\RP^m,N),\Theta \otimes \Q),$ $g \in S(N)$, and $\alpha \mapsto g(\alpha)$ is the induced action in the homology. Therefore the statement of Theorem \ref{main2}(B) concernong the case of odd $N$ follows from Theorem \ref{main1}.

The statement concerning even $N$ follows immediately from the exact homological sequence of the pair $(B(\RP^m,N), B(\RP^m_\star,N))$ with coefficients in the local system $ \tilde \Theta \otimes \pm \Q$, Theorem \ref{main2}A, and the
Thom isomorphism $H_i(B(\RP^m,N), B(\RP^m_\star,N); \tilde \Theta \otimes \pm \Q) \simeq H_{i-m}(B(\RP^m_\star, N-1), \tilde \Theta \otimes \pm \Q)$ for the normal bundle of the submanifold $B(\RP^m,N) \setminus B(\RP^m_\star,N) \simeq B(\RP^m_\star, N-1)$. \hfill $\Box$

\begin{lemma} 
\label{lll}
If $N$ is even and $m$ is odd, then the manifold $V(N) \subset B(\R^m,N) \subset B(\RP^m,N)$ participating in Lemmas \ref{lemW} and \ref{realiz} defines a non-trivial element of the group $H_{\frac{N}{2}(m-1)}(B(\RP^m,N), \tilde \Theta \otimes \pm \Q)$.
\end{lemma}

\noindent
{\it Proof} is very similar to that of Lemma \ref{realiz}: our cycle has non-zero intersection index with the oriented subvariety in $B(\RP^m,N)$ that consists of configurations which can be split into the pairs of points lying in the same fibers of the Hopf fibration. \hfill $\Box$

\section{On configuration spaces of $\RP^m$ and $\RP^m_\star$ with even $m$}

In addition to Theorems \ref{main2} and \ref{lem29}, we will use the following facts on the homology of unordered configuration spaces.

\begin{proposition}
\label{prepe}
If $m$ is even, then the groups $\bar H_i(I(\RP^m,N),\Q)$ and $\bar H_i(I(\RP_\star^m,N-1),\Q)$ are trivial for all $N \geq 2$ and $i\geq 0$. 
\end{proposition}

\noindent
{\it Proof.} As in \S \ref{sire}, consider the set $\nabla \equiv (\RP^m)^N \setminus I(\RP^m,N)$ and its simplicial resolution ${\mathfrak S} \subset \Delta(N) \times \nabla$, so that $\bar H_i(I(\RP^m,N),\Q) \simeq H_i((\RP^m)^N, \nabla; \Q)$ and $H_*(\nabla,\Q) \simeq H_*({\mathfrak S},\Q)$. Let us calculate the latter group by the spectral sequence built of the Borel--Moore homology groups of blocks (\ref{block}) with constant $\Q$-coefficients.
All spaces $L(A)$ have rational homology of a point, therefore this spectral sequence coincides with the one calculating the rational homology of the order complex $\Delta(N)$ of all non-complete partitions of $\{1, \dots, N\}$. This order complex has a maximal element (corresponding to the
``partition'' into a single set), hence is contractible. 
Therefore the space $\nabla$ also has the rational homology of a point, and our statement concerning the group $\bar H_*(I(\RP^m,N),\Q)$ follows from the exact sequence of the pair $\left((\RP^m)^N, \nabla \right)$. 

Similarly, for $N>1$ we have 
\begin{equation}
\label{exsec}
\bar H_i(I(\RP_\star^m,N-1),\Q) \simeq \tilde H_{i-1}(\tilde \nabla, \Q) \simeq \tilde H_{i-1}(\tilde {\mathfrak S}, \Q ) ,
\end{equation}
where $\tilde \nabla = \nabla \cap (\RP^m)^{N-1}$, see the previous section. The set $ \tilde {\mathfrak S} $ also has the rational homology of a point by a similar spectral sequence, so the groups (\ref{exsec}) are trivial. 
\hfill
$\Box$ 

\begin{corollary}
\label{refe}
If $m$ is even then the groups $\bar H_i(B(\RP^{m},N), \Q),$ $\bar H_i(B(\RP^{m},N), \pm \Q),$ 
$\bar H_i(B(\RP_\star^{m},N-1), \Q),$ and $\bar H_i(B(\RP_\star^{m},N-1), \pm \Q)$
are trivial for any $N \geq 2$ and $i \geq 0$. 
\end{corollary}

\noindent
{\it Proof}. These groups are the quotients of the groups considered in Proposition \ref{prepe} by some actions of permutation groups of $N$ or $N-1$ points. \hfill $\Box$ \medskip

\begin{lemma}
\label{lem9}
If $m$ is even then the group $\bar H_i(B(\RP^{m},N), \tilde \Theta \otimes \Q)$, $N\geq 2$, is equal to $\Q$ if $i=N m$ or $i=N m-(2m-1)$, and is trivial otherwise.
\end{lemma}

\noindent
{\it Proof}. By Remark \ref{protriv}, the orientation sheaf of $B(\RP^{m},N)$ is in this case equal to $\tilde \Theta$. Therefore our statement is Poincar\'e dual to the next one.

\begin{lemma}[see \cite{knu}]
\label{propl}
For any $N \geq 2$, the group $H_i\left(B(\RP^{m},N),\Q\right)$, $m$ even, is equal to $\Q$ if $i=0$ or $i=2m-1$, and is trivial for all other $i$. \hfill $\Box$
\end{lemma}

Let us describe the bases for the $\Q$-vector spaces considered in the last two lemmas. We consider $\RP^m$ as a quotient space of the $m$-sphere with the induced metric.

\begin{lemma}
\label{reali}
The vector space $H_{2m-1}\left(B(\RP^m,N),\Q\right)$ $($where $N\geq 2$ and $m$ is even$)$ is generated by the class of the submanifold in $B(\RP^m,N)$ consisting of configurations, all whose points lie in one and the same projective line in $\RP^m$ and separate this line into $N$ equal segments. The group $\bar H_{Nm-(2m-1)}(B(\RP^m,N), \tilde \Theta \otimes \Q)$ is generated by the class of the variety consisting of all configurations containing the fixed point $X_0 \in \RP^m$ and some other point lying in a fixed line passing through $X_0$.
\end{lemma}

\noindent
{\it Proof.} It is easy to see that these cycles are orientable in the first case and $\tilde \Theta$-orientable in the second one.
They have only one common point $C$: the first cycle is smooth at it, and the second one coincides in its neighbourhood with the union of $N-1$ smooth manifolds. Let us calculate the intersection index of these cycles at this point. 

Choose the chart $\R^m$ with coordinates $x_1, \dots, x_m$ in $\RP^m$ in such a way that the distinguished point $X_0$ is the origin in $\R^m$, the line containing the configuration $C$ is the $x_1$-axis, and remaining $N-1$ points of $C$ lie in this axis in the domain where $x_1>0$. Any line neighboring to this axis can be defined by the system of equations $x_k=p_k+q_kx_1$, $k=2, \dots, m$. For the local parameters of the first cycle in a neighborhood of the point $C$ we can choose the coefficients $p_k$ and $q_k$ of these equations of the line containing a configuration, and also the number $\varkappa$ defined as the coordinate $x_1$ of the point of this configuration neighboring to the origin. 

A regular point of the second cycle is a $N$-configuration containing the origin in $\R^m$ and one other point of the $x_1$-axis. The set of such regular configurations having a positive coordinate $x_1$ of this second point is path-connected. A coorientation of this cycle at such a point in $B(\R^m,N)$ is defined by the differential form \begin{equation}\label{coor}dx_1(u) \wedge dx_2(u) \wedge \dots \wedge dx_m(u) \wedge dx_2(v) \wedge \dots \wedge dx_m(v)\end{equation} where $u$ is the point of a neighboring configuration which is close to the origin, and $v$ is the other point of this configuration which is close to the $x_1$-axis. Moreover, this coorientation can be continued to the self-intersection points of the second cycle, corresponding to configurations with many points on the $x_1$-axis, as a collection of coorientations of all smooth local components of this cycle.

Let $0 \leq a_2 \leq \dots \leq a_N$ be the $x_1$-coordinates of points of the configuration $C$. Consider the square matrix of order $2m-1$, whose cell $(i,j)$ contains the value which the $i$-th factor of the coorientation form (\ref{coor}) of the local component of the second cycle related with the point $a_k$ takes on the $j$-th basic tangent vector of the first cycle at $C$ defined by increasing its $j$-th local parameter and keeping all other parameters. This matrix is equal to
$\mbox{
\begin{tabular}{|c|c|c|}
\hline
0 & 0 & 1 \\
\hline
 {\bf E} & {\bf 0} & 0\\
 \hline
 {\bf E} & $a_k${\bf E} & 0\\
\hline
\end{tabular}}$ \ 
where {\bf E} is the unit $(m-1)\times (m-1)$ matrix, {\bf 0} is zero matrix, and zeros in the top row and right-hand column mean strings or columns of $m-1$ zeros. In particular, its determinant is equal to $a_k^{m-1}$. All numbers $a_k$ are positive, therefore all these determinants are positive too, and the sum of their contributions to the intersection index of our two cycles is equal to $N-1 \neq 0$. In particular, the homology classes of these cycles are non-trivial. \hfill $\Box$ \medskip

Consider now the groups $H^*(B(\RP^m_\star,N), \Q)$ with even $m$. The space $\RP^m_\star$ is obviously fibered over the manifold $\RP^{m-1}$ with fiber $\R^1$ . Let us fix an orientation of $\RP^{m-1}$ and define for any $N$ the class $\mbox{Ind} \in H^{m-1}(B(\RP^m_\star,N))$, whose value on a generic piecewise-smooth $(m-1)$-dimensional cycle in $B(\RP^m_\star,N)$ is equal to the oriented number of $N$-configurations in this cycle such that their images in $\RP^{m-1}$ contain the distinguished point of $\RP^{m-1}$. It is easy to check that this cohomology class is well-defined.

\begin{lemma}
\label{porp}
The group $H^j(B(\RP^m_\star,N),\Q),$ $m$ even, is equal to $\Q$ for $j=0$ and $j=m-1$ and is trivial for all other $j$. The group $H^{m-1}(B(\RP^m_\star,N),\Q)$ is generated by the cocycle $\mbox{\rm Ind}$.
\end{lemma}

\noindent
{\it Proof} uses the induction over $N$. For $N=1$ the statement is obvious, since $\RP_\star^m$ is homotopy equivalent to $\RP^{m-1}$. The set $B(\RP^m,N) \setminus B(\RP^m_\star,N)$ can be obviously identified with the manifold $B(\RP^m_\star,N-1)$. All fibers of the normal bundle of this submanifold in $B(\RP^m,N)$ are canonically isomorphic to the tangent space of $\RP^m$ at the distinguished point, in particular this normal bundle is trivial. Therefore we have an isomorphism
\begin{equation} \label{ghys}
H_j(B(\RP^m,N), B(\RP^m_\star,N);\Q) \simeq H_{j-m}(B(\RP^m_\star,N-1),\Q). \end{equation}
Suppose now that our lemma is proved for the space $B(\RP_\star^m,N-1)$. By the exact sequence of the pair $(B(\RP^m,N), B(\RP^m_\star,N)),$ the groups $H_j(B(\RP_\star,N),\Q)$ have a chance to be non-trivial only for $j$ close to the dimensions of non-trivial elements of the groups $H_*(B(\RP^m,N),\Q)$ or $H_*(B(\RP^m,N),B(\RP^m_\star,N);\Q)$. 

By Lemma \ref{propl} and the induction hypothesis, there are only two fragments of this exact sequence, containing such non-trivial elements. One of them is 
$$
0 \to H_{2m-1}(B(\RP^m_\star,N),\Q) \to H_{2m-1}(B(\RP^m,N),\Q) \stackrel{\rho}{\longrightarrow} $$
$$ \stackrel{\rho}{\longrightarrow} H_{2m-1}(B(\RP^m,N), B(\RP^m_\star,N);\Q) \to H_{2m-2}(B(\RP_\star^m,N),\Q) \to 0 \ . $$
Two groups in the middle of this fragment, which are connected by the map $\rho$, are isomorphic to $\Q$ by Lemma \ref{propl}, isomorphism (\ref{ghys}) and the induction hypothesis. This map $\rho$ sends the class of the basic cycle described in the first statement of Lemma \ref{reali} to a class which corresponds via (\ref{ghys}) to the element of the group $H_{m-1}(B(\RP_\star^m,N-1),\Q),$ on which the cocycle $\mbox{Ind}$ takes the value of $\pm (N-1) \neq 0$. So the map $\rho$ is an isomorphism, and all neighboring terms of our exact sequence vanish.

Another suspicious fragment of our sequence is $$0 \to H_m(B(\RP^m,N), B(\RP^m_\star,N);\Q) \to H_{m-1}(B(\RP^m_\star,N), \Q) \to 0 \ . $$ Its group $H_m(B(\RP^m,N), B(\RP^m_\star,N);\Q)$ is isomorphic by (\ref{ghys}) to the group $H_0(B(\RP_\star^m,N-1),\Q) \simeq \Q$. Hence the group $ H_{m-1}(B(\RP^m_\star,N), \Q)$ is also isomorphic to $\Q$. \hfill $\Box$

\section{Proof of Theorems \ref{cor1} and \ref{m4}} 

We know now all groups in the formulas (\ref{ssss}), (\ref{thomm2}) and (\ref{sssspt}--\ref{lenspmpt}) tensored with $\Q$ for all combinations of parities of $m$ and $M$. Let us apply this information. 

\subsection{Cases of even $m$ and even $i(M+1)$}
If $M$ is odd then the coefficient sheaves in both (\ref{sssspt}) and (\ref{thomm2pt}) are equal to $\pm \Z$. If $M$ is even then
the coefficient sheaf in (\ref{sssspt}) is equal to $\Z$.
By Corollary \ref{refe} (see page \pageref{refe}), in all these cases $\hat E_1^{p,q} \otimes \Q \equiv 0$ for all $p < 0$. 
This implies that all spaces $\mbox{EM}^*_0(S^m,S^M)$ with even $m$, and also all spaces $\mbox{EM}^*_1(S^m,S^M)$ with even $m$ and odd $M$ have the rational homology of a point. By the spectral sequence of the fiber bundle (\ref{mfb}), the Poincare polynomials of corresponding spaces $\mbox{EM}_0(S^m,S^M)$ and $\mbox{EM}_1(S^m,S^M)$ are then equal to $1+ t^M$.

\subsection{The case of even $m$, even $M$, and $i=1$}
\label{ee1}

The coefficient sheaf in (\ref{thomm2}) is in this case equal to $\tilde \Theta$.
By Lemma \ref{lem9}, there are then only the following non-trivial groups $E^{p,q}_1 \otimes \Q$ with $p<0$ (all of which are isomorphic to $\Q$): they correspond to $(p,q) = (-1,M-m+1)$, or $p\leq -2$ and either $q=-p(M-m+1)$ or $q=-p(M-m+1)+(2m-1)$. 
Let us prove that $E_\infty \equiv E_1$ over $\Q$. 

The group $E^{-1,M-m+1} \simeq \Q$ obviously survives until $E_\infty$ and provides a generator $[\Sigma]$ of the group $H^{M-m}(\mbox{EM}_1(S^m,(\R^{M+1} \setminus \{0\}),\Q )\simeq \Q$, namely the linking number class of the entire {\em discriminant variety} in the space of all odd maps $S^m \to \R^{M+1}$. (This variety consists of maps whose images contain the point $0 \in \R^{M+1}$; its codimension is equal to $M-m+1$). It is enough to prove that all powers $[\Sigma]^{\smile N}$ of this basic class $[\Sigma]$ are non-trivial elements of the groups $H^{N(M-m)}(\mbox{EM}_1(S^m,(\R^{M+1} \setminus \{0\})) ,\Q )$. Consider a generic point of the $N$-fold self-intersection of the discriminant variety, that is, some odd $C^\infty$-map $F:S^{m} \to \R^{M+1}$ such that the set $F^{-1} (0)$ consists of exactly $N$ pairs of opposite points, and $dF$ is injective at all these points. Let ${\mathcal L} \ni F$ be a generic affine $(N (M+1))$-dimensional subspace in the space of $C^\infty$-smooth odd maps $S^{m}\to \R^{M+1}$. Close to $F$, the intersection of ${\mathcal L}$ with the discriminant variety consists of $N$ locally irreducible smooth components meeting generically one another. The complement of this intersection in a small neighborhood of the point $F$ in ${\mathcal L}$ is homotopy equivalent to the product of $N$ spheres of dimension $M-m$. Any factor of this product has linking number 1 with only one local component of the discriminant variety $\Sigma \cap {\mathcal L}.$ So the restriction of the globally defined class $[\Sigma] \in H^{M-m}(\mbox{EM}_1(S^{m},(\R^{M+1} \setminus 0)) ,\Q )$ to any of these factors takes value 1 on the (suitably oriented) fundamental class of this factor. Since the number $M-m$ is even, the $N$-th power of the class $[\Sigma]$ takes value $N!$ on the fundamental $(N(M-m))$-dimensional cycle of this product. In particular, it is not zero-cohomologous. Therefore the additional cells $E^{p,q}$ with $p<-1$, $q=-p(M-m+1)+(2m-1)$, also survive in $E_\infty$.
The desired Poincar\'e series is thus equal to 
\begin{equation}
1+ t^{M-m} + \sum_{N \geq 2} t^{N(M-m)} (1+t^{2m-1}) \equiv \frac{1+ t^{2M-1}}{1-t^{M-m}} \ . 
\label{3evev}
\end{equation}

In a similar way, Lemma \ref{porp} and formula (\ref{thomm2pt}) imply that all non-trivial groups $\hat E^{p,q}_1$ of our spectral sequence calculating the group $H^*(\mbox{EM}^*_1(S^m,S^M),\Q)$ are equal to $\Q$ and correspond to the pairs $(p,q)$ equal to $(0,0)$ and $(-N,N(M-m+1))$ and $(-N,N(M-m+1)+m-1)$ for arbitrary $N \geq 1$. The stabilization $\hat E_\infty \equiv \hat E_1$ of this spectral sequence can be proved exactly as in the previous paragraph. The Poincare polynomial of the resulting homology group is thus equal to
\begin{equation}
1+ \sum_{N \geq 1} t^{N(M-m)} (1+t^{m-1}) \equiv \frac{1+ t^{M-1}}{1-t^{M-m}} \ . \label{2evev}
\end{equation} 
The formulas (\ref{3evev}) and (\ref{2evev}) form the third and the second cells of the corresponding row of Table \ref{oddd}. The structure of the differential $d^M$ indicated in the last cell of this row is the unique one compatible with these formulas.

\subsection{The case of odd $m$, odd $M$ and arbitrary $i$}
 
By Lemma \ref{orient} the orientation sheaf of $B(L_r^{m},N)$ is equal to $\pm \Z$. This is exactly the coefficient sheaf in (\ref{lenspm}). Hence the group (\ref{lenspm}) is Poincar\'e isomorphic to the group 
$ H^{q+p(M-m+1)}(B(L_r^{m},-p), \Z).$
By Theorem \ref{lem29}(a) all non-trivial groups $E_1^{p,q} \otimes \Q,$ $p<0$, are in this case the ones with arbitrary $p \leq -1$ and either $q= -p(M-m+1)$ or $q=-p(M-m+1) +m$; all of them are equal to $\Q$. All these groups survive until $E_\infty$ for the same reasons as in the previous case.
So, the Poincar\'e series of the group $H_*(\mbox{EM}_i(S^m,S^M),\Q)$ with arbitrary $i$ is equal to 
$$1+ \sum_{N \geq 1} t^{N(M-m)} (1+t^{m}) \equiv \frac{1+ t^{M}}{1-t^{M-m}} \ , $$ 
see the top cells of the third columns in Tables \ref{evv} and \ref{oddd}, and the first statement of Theorem \ref{m4}.

Also, the group (\ref{lenspmpt}) is Poincare isomorphic to $H^{q+p(M-m+1)}(B(L^m_{r,\star} \ ,-p),\Z)$. 
By statement (\ref{b}) of Lemma \ref{lem11} and Proposition \ref{chn76}, the latter group tensored with $\Q$ is equal to $\Q$ if $q+p(M-m+1)=0$ and is trivial in all other cases. Therefore we have $\hat E^{p,q}_1 \otimes \Q \simeq \Q$ for $(p,q) = (-N, N(M-m+1))$, $N \geq 0$; all other groups $\hat E^{p,q}_1 \otimes \Q$ with $p<0$ are trivial. This gives us the Poincare series $$\sum_{N=0}^\infty t^{N(M-m)} \equiv \frac{1}{1-t^{M-m}} \ ,$$ see the top cells of the second columns in Tables \ref{evv} and \ref{oddd}, and also the second statement of Theorem \ref{m4}. 

\subsection{The case of odd $m$, even $M$, and $i=0$}
\label{3eoe}
The corresponding group $E^{p,q}_1$ given by (\ref{ssss}) is Poincar\'e isomorphic to
$ H^{q+p(M-m+1)}(B(\RP^{m},-p), \pm \Z).$

By Theorem \ref{lem29}(B), the non-trivial groups $E^{p,q}_1 \otimes \Q$, $p<0$, lie on two parallel lines in the $(p,q)$-plane. Namely,
$E^{p,q}_1 \otimes \Q \simeq \Q$ if $p$ is odd and the number $q+p(M-m+1)$ is equal to either $(-p-1)(m-1)/2$ or to $(-p-1)(m-1)/2+m$; in all other cases $E^{p,q}_1 \otimes \Q =0$. Since $M>m$, no differential in this sequence can connect two non-trivial groups, hence $E_\infty \equiv E_1$ over $\Q$. 
By counting the total dimensions $p+q$ of these group for all values $p=-(2j+1),$ $j \geq 0$, we get the Poincar\'e series 
$$ 1+ \sum_{j=0}^\infty (1+t^m)t^{M-m +j(2M-m-1)} \equiv 1+ t^{M-m}\frac{1+t^m}{1-t^{2M-m-1}} \ , $$
see the third cell of the last row in Table \ref{evv}.

By (\ref{sssspt}) and Poincar\'e duality, $\hat E^{p,q}_1 \otimes \Q \simeq H^{q+p(M-m+1)}(B(\RP_\star^{m},-p), \pm \Q).$ 
By statement (\ref{c}) of Lemma \ref{lem11} and Proposition \ref{chn76}, this group is equal to $\Q$ if $q+p(M-m+1) = (m-1)[\frac{N}{2}]$ and is trivial in all other cases. So for any $N>0$ there is exactly one non-trivial cell $E^{-N,q}_1$, namely $q$ should be equal to $N(M-m+1)+(m-1)[N/2]$. This expression depends monotonically on $N$, therefore our spectral sequence has no non-trivial higher differentials, and $\hat E_\infty = \hat E_1$. The Poincare series of the resulting homology group is thus equal to 
$$1 + t^{M-m} + \sum_{k=1}^\infty \left(t^{2k(M-m)+k(m-1)} + t^{(2k+1)(M-m)+k(m-1)}\right) \equiv \frac{1+t^{M-m}}{1-t^{2M-m-1}} \ ,$$
see the second cell of the last row in Table \ref{evv}. 

The structure of differentials $d^M$ indicated in the last cell of this row is the unique one which is compatible with the previous two cells.

 \subsection{The case of odd $m$, even $M$, and $i=1$} 
\label{3eoo}
Substituting Theorem \ref{main2}(A) to (\ref{thomm2}), we obtain that the group $E_1^{p,q} \otimes \Q$ is trivial if $p$ is odd; if $p=-N \neq 0$ is even, then $E_1^{p,q} \otimes \Q$ is equal to $\Q$ for $q=N(M+1)- \frac{N}{2}(m+1)$ or $q=N(M+1)- \frac{N}{2}(m+1)+m$, and is trivial for all other values of $q$. Again, our rational spectral sequence stabilizes at the term $E_1$ by dimensional reasons. The counting of the total degrees $p+q\equiv -N+q$ of its groups over all even natural values of $N$ gives us the Poincar\'e series equal to $$1 + (1+t^m)\sum_{j=1}^\infty t^{j(2M-m-1)} \equiv \frac{1+t^{2M-1}}{1-t^{2M-m-1}}.$$ 

Substituting Theorem \ref{main2}(B) to (\ref{thomm2pt}), we obtain that the group $\hat E_1^{p,q} \otimes \Q$ is non-trivial (and then equal to $\Q$) only for $(p,q)=\left(-N,N(M-m+1)+\left]\frac{N}{2}\right[\times (m-1)\right)$, $N \geq 0$. The stabilization $\hat E_\infty \equiv \hat E_1$ follows immediately from the position of these groups. Counting the corresponding total degrees $p+q$ separately forr even and odd $N$ we get the Poincar\'e series $$ \sum_{k=0}^\infty t^{k(2M-m-1)} + \sum_{k=0}^\infty t^{k(2M-m-1)+M-1} \equiv \frac{1+t^{M-1}}{1-t^{2M-m-1}} .$$ So we have found the second and the third cells of the last row of Table \ref{oddd}. The structure of the differential $d^M$ indicated in the last cell of this row is the only one compatible with the previous two cells. 

\section{The justification of Table \ref{genn}}

By Poincare duality and Lemma \ref{orient}, the group (\ref{msgen}) is equal to 
\begin{equation}
\label{msgen1}
H^{q+p(M-m+1)}\left(B(S^m,-p),(\pm \Z)^{\otimes(M-m)}\right) ,
\end{equation}
and the group (\ref{msgenpt}) to 
\begin{equation}
\label{msgen1pt}
H^{q+p(M-m+1)}\left(B(\R^m,-p),(\pm \Z)^{\otimes(M-m)}\right) .
\end{equation}

\subsection{$m$ and $M$ are odd}

See Theorem \ref{m4}, the case of $r=s=1$.

\subsection{$m$ and $M$ are even}

The group $H^*(B(\R^m,N),\Q)$ is in this case equal to $\Q[0]$ for $N=1$ and to $\Q[0] \oplus \Q[m-1]$ for $N \geq 2$ (see \cite{Cohen76}, Section 4). Consequently, all non-trivial cells $\hat E^{p,q}_1 \otimes \Q$, $p <0,$ of the spectral sequence $\hat E_r^{p,q}$ correspond to the pairs $(p,q)$ equal to either $(-1, M-m+1)$ or to $(-N,N(M-m+1))$ and $(-N, N(M-m+1)+m-1)$ with arbitrary $N \geq 2$; all these non-trivial groups are isomorphic to $\Q$. All of them survive until $\hat E_\infty$ by the same reason as in \S \ref{ee1}. So the Poincare series of the resulting homology group is equal to $$1+ t^{M-m} + \sum_{N=2}^\infty \left(t^{N(M-m)} + t^{N(M-m)+m-1}\right) \equiv \frac{1+t^{2M-m-1}}{1-t^{M-m}} \ ,$$ see the second column of Table \ref{genn}. 

The group $H^*(B(S^m,N),\Q)$ is equal to $\Q[0] \oplus \Q[2m-1]$ for $m \geq 3$
by the exact homological sequence of the pair $(B(S^m,N), B(\R^m,N))$ (where $\R^m$ is $S^m$ less one point).
Indeed, the difference $B(S^m,N) \setminus B(\R^m,N)$ 
is obviously diffeomorphic to $B(\R^m,N-1)$, and the normal bundle of this difference in $B(S^m,N)$ is constant, with the fiber equal to the tangent space of $S^m$ at the removed point. Therefore by the Thom isomorphism and the well-known structure of $H^*( B(\R^m,N-1),\Q)$ the group $H_i(B(S^m,N), B(\R^m,N);\Q)$ is trivial in all dimensions except for $m$ and $2m-1$ and is isomorphic to $\Q$ in these two dimensions. The boundary map $H_i(B(S^m,N), B(\R^m,N),\Q) \to H_{i-1}(B(\R^m,N),\Q)$ is thus trivial for all $i \neq m$. If $i=m$, then by the construction of the Thom isomorphism the image of this map is generated by the $(m-1)$-dimensional submanifold in $B(\R^m,N)$ swept out by configurations, some $N-1$ points of which are constant, and the $N$-th one runs a huge sphere containing these $N-1$ points inside. This cycle is non-trivial in $H_{m-1}(B(\R^m,N),\Q)$, so our boundary map is also non-trivial for $i=m$.

Finally, we get the following structure of the term $E_1$ of the main spectral sequence over $\Q$. All groups $E_1^{p,q}$ with the following pairs $(p,q)$ are equal to $\Q$: $(0,0)$ (the standard term), $(-1,M-m+1)$ and $(-1,M+1)$ (because $B(S^m,1) = S^m)$), $(-2,2(M-m+1)$ (because $B(S^m,2)$ is homotopy equivalent to $\RP^m)$), and $(-N,N(M-m+1))$ and $(-N,N(M-m+1)+2m-1)$ for any $N \geq 3$. 
All other groups $E_1^{p,q} \otimes \Q$ are trivial. The spectral sequence stabilizes at this term, that is $E_\infty \equiv E_1$, by the same reasons as in \S \ref{ee1}. Therefore the Poincare series of the group $H_*(\mbox{Map}(S^m,S^M),\Q)$ is equal to
$$1 + t^{M-m} + t^M + t^{2(M-m)} + \sum_{N=3}^\infty \left(t^{N(M-m)} + t^{N(M-m)+2m-1}\right),$$ which is equal to the expression in the third column of Table \ref{genn}. The form of the higher differentials indicated in the fourth column for the case of even $m$ and $M$ is the unique one compatible with the data of the second and third cells.

\subsection{$m$ is even, $M$ is odd}

If $m$ is even then the group $H^*(B(\R^m,N),\pm \Q)$ is trivial for any $N \geq 2$, as well as its Poincar\'e dual group $\bar H_*(B(\R^m,N),\pm \Q)$, see e.g. Corollary 2 in \S I.4 of \cite{Book}. Therefore the spectral sequence (\ref{msgenpt}) for such $m$ and $M$ has only one non-trivial group $E^{p,q}_1$, with $p<0$, namely $E^{-1,M-m+1}_1$. This gives us the Poincare series $1+t^{M-m}$ for the group $H_*(\Omega^m S^M,\Q)$. The spectral sequence of the fiber bundle (\ref{fb}) has therefore only four non-trivial cells: ${\mathcal E}_2^{p,q} \simeq \Q$ for $p \in \{0,M\}$ and $q \in \{0,M-m\}$. It obviously stabilizes in this term ${\mathcal E}_2$ and gives us the corresponding cell of the third column of Table \ref{genn}.

\subsection{$m$ is odd, $M$ is even}

By Proposition \ref{chn76} the group $H^j(B(\R^m,N), \pm \Q) \equiv \bar H_{mN-j}(B(\R^m,N), \Q)$ is in this case non-trivial (and then equal to $\Q$) only for $j=(m-1)[N/2].$ So, for any $N>0$ there is exactly one non-trivial cell $\hat E^{-N,q}_1$, namely the one with $q=N(M-m+1)+(m-1)[N/2]$. The number $q$ in this expression depends strictly monotonically on $N$, therefore all higher differentials are trivial, and $\hat E_\infty = \hat E_1$. By counting the total degrees $p+q$ of non-trivial groups $\hat E_\infty^{p,q}$ separately for even and odd $N$ we see that the Poincar\'e series of the group $H_*(\Omega^m S^M,\Q)$ is equal to $$\sum_{k=0}^\infty \left(t^{2k(M-m)+k(m-1)} + t^{(2k+1)(M-m)+k(m-1)}\right) \equiv \frac{1+t^{M-m}}{1-t^{2M-m-1}},$$
see the second cell of the last row in Table \ref{genn}. 

Further, by Theorem \ref{lem29}(B) (with $r=1$) all groups $H^j(B(S^m,N), \pm \Q)$ are trivial if $N$ is even; if $N$ is odd then they are non-trivial (and then equal to $\Q$) only for $j$ equal to $(m-1)\frac{N-1}{2}$ or $(m-1)\frac{N-1}{2}+m$. Therefore the calculation of the main spectral sequence $E_r^{p,q}$ in this case coincides with the one given in the first paragraph of \S \ref{3eoe}.

\section{The construction of the spectral sequences of Theorems \ref{absv} and \ref{relav}}
\label{sseven}

This construction follows very closely the strategy described in Chapter 3 of \cite{Book}. Namely, we consider the {\em discriminant} set $\Sigma \equiv \mbox{EM}_G(X,\R^W) \setminus \mbox{EM}_G(X,(\R^W \setminus C\Lambda))$, i.e. the set of $G$-equivariant maps $X \to \R^W$ whose images meet $C\Lambda$. For any finite-dimensional affine subspace ${\mathcal F} \subset \mbox{EM}_G(X,\R^W)$ we study the group $\bar H_*( \Sigma \cap {\mathcal F})$, which is Alexander dual to the group $H^*({\mathcal F} \setminus \Sigma)$, cf. \cite{A70}. The groups $H^*({\mathcal F}_i \setminus \Sigma)$ of an increasing sequence of generic subspaces ${\mathcal F}_i$ stabilize to the homology group of the space $\mbox{EM}_G(X,(\R^W \setminus C\Lambda))$ by a variant of the Weierstrass approximation theorem. On the other hand, the groups $\bar H_*({\mathcal F}_i \cap \Sigma)$ can be calculated using spectral sequences which stabilize (after a change of indices reflecting the Alexander duality) to the spectral sequence of Theorem \ref{absv}. 

Everywhere in this section we use the notation and assumptions of \S \ref{agen}. In particular, $X$ is a compact semialgebraic subset of $\R^a$, $G$ is a compact Lie group acting freely and algebraically on $X$ and on some neighborhood of $X$ in $\R^a$, $\rho: G \to O(\R^W)$ is a representation of $G$, and $C\Lambda$ is a $\rho(G)$-invariant cone in $\R^W$.

\subsection{Finite-dimensional approximations}
\label{fda}

\begin{proposition}
\label{appro}
For any homology class $\alpha \in H_*( \mbox{\rm EM}_G(X,(\R^W \setminus C\Lambda))),$ there is a finite-dimensional affine subspace $\mathcal{F} \subset \mbox{\rm EM}_G(X,\R^W)$ such that the class $\alpha$ can be realized by a cycle contained in the subspace ${\mathcal F} \setminus \Sigma \equiv {\mathcal F} \cap \mbox{\rm EM}_G(X,(\R^W \setminus C\Lambda)). $ 

If two cycles in such a subspace ${\mathcal F} \setminus \Sigma$ represent one and the same element of the group $H_*( \mbox{\rm EM}_G(X,(\R^W \setminus C\Lambda))),$ then there is a finite-dimensional subspace ${\mathcal F}' \supset {\mathcal F}$ in $\mbox{\rm EM}_G(X,\R^W)$
such that these two cycles are homologous to one another in ${\mathcal F}' \setminus \Sigma$. 

The relative versions of these statements for the spaces \ $\mbox{\rm EM}_G(X,A;(\R^W \setminus C\Lambda))$ are also true.
\end{proposition}

\noindent
{\it Proof.} Realize the class $\alpha$ by a continuous map $\varkappa: \Xi \to \mbox{EM}_G(X,(\R^W \setminus C\Lambda))$ of a compact polyhedron $\Xi$ (i.e. by a continuous map $\tilde \varkappa: X \times \Xi \to \R^W \setminus C\Lambda$ that is equivariant with respect to the $G$-actions on $X$ and $\R^W$ and such that $\tilde \varkappa(x , \xi) \equiv \varkappa(\xi)(x)$ for $x \in X, \xi \in \Xi$). Let $r$ be the minimal distance in $\R^W$ between the sets $f(X)$ and $C\Lambda$ over all points $f \in \varkappa(\Xi)$. Let us embed the complex $\Xi$ regularly into some Euclidean space $\R^b$. $\frac{r}{2}$-approximate the map $\tilde \varkappa$ by the restriction of a polynomial map $\R^a \times \R^b \to \R^W$. This approximation can be considered as an $\frac{r}{2}$-approximation (in the $C^0$-topology) of the cycle $\varkappa(\Xi)$ by a cycle in the finite-dimensional vector space of restrictions to $X$ of polynomial maps $F: \R^a \to \R^W$ of a certain finite degree. Then replace any such map $F|_X$ by its $G$-symmetrization (that is, by the map $X \to \R^W$ sending any point $x \in X$ to $$ \int_G\rho(g^{-1})(F(g(x))) d \mu \ , $$
integration over the Haar measure). This $G$-symmetrization is a linear operator from the finite-dimensional space of polynomial maps of a certain degree to a subspace of the space of $G$-equivariant maps. The latter subspace is thus also finite-dimensional. The image in it of the cycle $\varkappa(\Xi)$ after all perturbations $\frac{r}{2}$-approximates the initial cycle, in particular is homologous to it in $\mbox{EM}_G(X,(\R^W \setminus C\Lambda))$.

The second statement of Proposition \ref{appro} can be proved in a similar way by approximating homologies between our cycles in the space of equivariant maps.

In the relative case of Theorem \ref{relav} we use the following modification of this construction. Lift our inclusion $ X \hookrightarrow \R^a$ to the embedding $X \to \R^{a} \oplus \R^1$ sending any point $x \in X$ to $(x, \mbox{dist}(x,A))$.
Let $\Phi:X \to \R^W$ be an arbitrary $G$-equivariant extension of the given map $\varphi: A \to \R^W$ to entire $X$.
Again realize the homology class $\alpha$ by a map $\varkappa:\Xi \to \mbox{EM}_G(X,(\R^W \setminus C\Lambda)),$ then consider the map $X \times \Xi \to \R^W$ that takes any pair $(x,\xi)$ to $\varkappa(\xi)(x) - \Phi(x),$ and approximate this map by a polynomial map $\R^a \oplus \R^1 \oplus \R^b \to \R^W$ given by a set of $W$ polynomials vanishing on the hyperplane $\R^a \times \{0\} \times \R^b$. We obtain an approximation of the cycle $\alpha$ by a cycle in the affine space of maps of the form $P + \Phi$, where $P$ are polynomials of uniformly bounded degree. Finally, $G$-symmetrize these maps $P+\Phi$ as previously.
 \hfill $\Box$

\subsection{$k$-sufficient finite-dimensional approximations}

\begin{definition} \rm
\label{kdef}
A finite-dimensional affine subspace ${\mathcal F} \subset \mbox{EM} _G(X,A; \R^W)$ is {\it $k$-sufficient} if

1) for any natural $N \leq k$, any points $x_1, \dots, x_N$ of $X\setminus A$ from pairwise different orbits of $G$, and any points $\lambda_1, \dots, \lambda_N$ in $C\Lambda$, the plane in ${\mathcal F}$ that consists of maps $f$ such that $f(x_j)=\lambda_j$, $j=1, \dots, N$, has codimension exactly $NW$ in ${\mathcal F}$; 

2) for any natural $N$ the set of maps $f \in {\mathcal F}$ such that $f(x) \in C\Lambda$ for some $N$ points $x \in X \setminus A$ from pairwise different orbits of $G$ has codimension $\geq N(W-\dim X +\dim G- \dim C\Lambda)$ in ${\mathcal F}$. 

A subspace ${\mathcal F}$ is called {\it $k$-transversal} if it satisfies condition 1) of this definition.
\end{definition}

We say that a subset of the space of $L$-dimensional subspaces of $\mbox{EM} _G(X,A; \R^W)$ is semialgebraic if its intersection with the space of $L$-dimensional subspaces of any finite-dimensional affine subspace of $\mbox{EM} _G(X,A; \R^W)$ is semialgebraic. Such a subset is {\em of codimension at least $T$} if for any finite-dimensional affine subspace $\Phi \subset \mbox{EM} _G(X,A; \R^W)$ there is a subspace $\tilde \Phi \supset \Phi$ such that the intersection of this subset with the space of $L$-dimensional affine subspaces of $\tilde \Phi$ has codimension at least $T$ in this space. It is easy to see that the complement of any semialgebraic subset of positive codimension is dense in the space of $L$-dimensional subspaces of $\mbox{EM} _G(X,A; \R^W)$, and the complement of any semialgebraic subset of codimension $\geq 2$ is path-connected.

\begin{proposition}
\label{ksuff}
Let $X$, $A$, $G$ and $\Lambda$ be as above. Then 

1. For any natural $k$ and $L\geq k(W+\dim X - \dim G+ \dim C\Lambda),$ the set of not $k$-transversal subspaces is a semialgebraic subset of codimension at least $L- k(W+\dim X - \dim G+ \dim C\Lambda)+1$ in the space of $L$-dimensional affine subspaces of $\mbox{\rm EM} _G(X,A; \R^W)$.

2. If $\dim X - \dim G + \dim C\Lambda < W$ then for any natural $L$ there is a dense subset in the space of affine $L$-dimensional subspaces of $\mbox{\rm EM} _G(X,A; \R^W)$ such that condition 2$)$ of Definition \ref{kdef} is satisfied for all elements of this subset. 

3. $k$-sufficient subspaces are dense in the space of affine subspaces of $\mbox{\rm EM} _G(X,A; \R^W)$ of any sufficiently large finite dimension. 
\end{proposition}

\noindent
{\it Proof.} 1. For any collection of points $x_1, \dots, x_N$ and $\lambda_1, \dots, \lambda_N$ as in condition 1) of Definition \ref{kdef}, the maps $f$ such that $f(x_i)=\lambda_i$ for all $i=1, \dots, N$ form a subspace of codimension $N\, W$ in $\mbox{\rm EM} _G(X,A; \R^W)$. This subspace depends only on the class of the collection $\{(x_i, \lambda_i)\}$ in the space $\Pi_{N,G}(X,A)$, so we have an $N(\dim X-\dim G + \dim C\Lambda)$-parametric family of such subspaces. The set of $L$-dimensional affine subspaces in $\mbox{\rm EM} _G(X,A; \R^W)$ that are not in general position with respect to some fixed subspace of this family (i.e. have an empty or non-transversal intersection with it) has codimension $L-NW+1$ in the space of all $L$-dimensional subspaces. The union of all these sets over all points of $\Pi_{N,G}(X,A)$  is semialgebraic by Tarski-Seidenberg theorem, and the codimension of this union at least $L-NW+1-N(\dim X -\dim G + \dim C\Lambda) \geq L-k(W+\dim X -\dim G + \dim C\Lambda)+1 $. The codimension of the union of these unions over all $N \leq k$ also satisfies this estimate. 

2. Set $k_0 = [L/(W-\dim X +\dim G - \dim C\Lambda)]+1$. Given an $L$-dimensional affine subspace ${\mathcal F}$ of $\mbox{\rm EM} _G(X,A; \R^W)$ and an arbitrary vicinity of ${\mathcal F}$ in the space of all subspaces of this dimension, we first choose a $k_0$-transversal affine subspace 
 $\tilde {\mathcal F} \subset \mbox{\rm EM} _G(X,A; \R^W)$ of finite dimension greater than $k_0(W+\dim X -\dim G +\dim C\Lambda)$, which contains some $L$-dimensional subspace ${\mathcal F}'$ from this vicinity: such a subspace $\tilde {\mathcal F}$ exists by statement 1 of our proposition. Then for any $N \leq k_0$ the 
set of maps $f \in \tilde {\mathcal F}$ such that $f(x_i) \in C\Lambda$ for some $N$ points $x_i \in X$ from pairwise distinct $G$-orbits is swept out by a $N(\dim X -\dim G+\dim C\Lambda)$-parametric family of subspaces of codimension $NW$, in particular it is a semialgebraic subvariety of codimension $\geq N(W-\dim X +\dim G - \dim C\Lambda)$ in $\tilde {\mathcal F}$. 
Any $L$-dimensional affine subspace ${\mathcal F}''$ of $\tilde {\mathcal F}$ that is transversal to these subvarieties for all $N \leq k_0$ satisfies condition 2) of Definition \ref{kdef}. Indeed, for all $N\leq k_0$ this condition follows from the transversality; for $N =k_0$ it means that the intersection of ${\mathcal F}''$ with the corresponding subvariety is empty, which implies the same also for all $N> k_0$. 

By transversality theorem, we can choose this subspace ${\mathcal F}''$ to be arbitrarily close to ${\mathcal F}'$ and hence also to ${\mathcal F}$.

3. If $L \geq k(W+\dim X-\dim G + \dim C\Lambda)$ then we can choose this subspace ${\mathcal F}''$ to be also $k$-transversal.
\hfill $\Box$ 

\subsection{The simplicial resolution}
\label{sre}

We will show now that the groups $H^*({\mathcal F} \setminus \Sigma)$ for all $k$-sufficient subspaces ${\mathcal F}$ in
$\mbox{EM}_G(X, \R^W)$ or $ \mbox{EM}_G(X,A; \R^W)$
can be approximated using the standard spectral sequences of Theorems \ref{absv} and \ref{relav} in
all dimensions strictly smaller than $ (k+1)(W-\dim X + \dim G - \dim C\Lambda-1)$. 

Let ${\mathcal F}$ be a $k$-sufficient subspace in $\mbox{EM}_G(X,\R^W)$. A simplicial resolution of the variety $\Sigma \cap {\mathcal F}$ is then constructed as follows.

Let $T$ be the maximal number of different orbits of $G$ in $X$ which are mapped to the points of $C\Lambda$ by a map of the class ${\mathcal F}$. By the second property of $k$-sufficient spaces, this number does not exceed \ $\dim {\mathcal F}/(W -\dim X +\dim G - \dim C\Lambda)$. 

The simplicial resolution will be constructed as a subset in the product of the space ${\mathcal F}$ and the $T$-th {\em self-join} of the quotient space $X/G$.

Namely, let us embed this quotient space $X/G$ generically into an affine space of a very large dimension, and for any $N \leq T$ points of the image of the embedding consider their convex hull. If the dimension of the ambient space is sufficiently large and the embedding is indeed generic, then any such convex hull is a simplex of the corresponding dimension $N-1$, and these simplices have no unexpected intersections (that is, the intersection of any two such simplices is their common face, which is the convex hull of the intersection of their vertex sets). Let us assume that these conditions are satisfied, and define the $T$-th self-join $(X/G)^{\star T}$ of $X/G$ as the union of all these simplices with arbitrary $N \leq T$. Obviously, such spaces $(X/G)^{\star T}$ defined using different generic embeddings are homeomorphic to one another. 
The space $(X/G)^{*T}$ is compact since $X/G$ is.

Let $f \in \Sigma \cap {\mathcal F}$ be a $G$-equivariant map $f: X \to \R^W$ that takes exactly $N$ orbits of our action of $G$ on $X$ to $C\Lambda$. Associate with it the simplex in the self-join $(X/G)^{\star T}$ spanned by the corresponding $N$ points of $X/G$, and the simplex in $ (X/G)^{\star T}\times {\mathcal F}$ equal to the product of the previous simplex and the point $\{f\} \in {\mathcal F}$. Finally, define the simplicial resolution $\sigma \subset (X/G)^{\star T} \times {\mathcal F} $ as the union of these simplices for all discriminant maps $f \in \Sigma \cap {\mathcal F}$.

This space $\sigma$ has a natural filtration $\sigma_1 \subset \dots \subset \sigma_{T} \equiv \sigma$, namely, $\sigma_N$ is the union of all $\leq N$-vertex simplices in this construction.

\begin{proposition}
\label{pro6}
1. The map $\sigma \to \Sigma \cap {\mathcal F}$ induced by the obvious projection $ (X/G)^{\star T} \times {\mathcal F} \to {\mathcal F}$ is a proper map that induces an isomorphism $\bar H_*(\sigma) \to \bar H_*(\Sigma \cap {\mathcal F})$ of Borel--Moore homology groups of these spaces.

2. The dimension of any non-empty set $\sigma_N \setminus \sigma_{N-1}$ does not exceed \ $\dim{\mathcal F}-N(W-\dim X + \dim G - \dim C\Lambda-1)-1$.

3. If $N \leq k$ then $\sigma_N \setminus \sigma_{N-1}$ is the space of a fibered product of two fiber bundles with base $\Pi_{N,G}(X)$ $($see \S \ref{agen}$)$. The fibers of these two bundles are an $(N-1)$-dimensional open simplex and a $(\dim {\mathcal F} -N\, W)$-dimensional affine space respectively. The orientation sheaves of these two bundles are respectively the sign sheaf $\pm \Z$ of the unordered configuration space $B(X/G,N)$, and the sheaf ${\mathcal J} \otimes \pm \Z$, see \S \ref{agen}.
\end{proposition}

\noindent
{\it Proof.} Properness of the projection $\sigma \to \Sigma \cap {\mathcal F}$   follows immediately from the construction, and the isomorphism in Statement 1 is a standard property of simplicial resolutions, see e.g. \cite{Book}, \cite{gorinov}. Statement 2 follows  from the second condition in the definition of $k$-sufficient subspaces.

The projection of the fiber bundle from statement 3 is defined as follows. Any point of the space $\sigma_N \setminus \sigma_{N-1}$ has the form $(s,f)$, where $s$ belongs to an $(N-1)$-dimensional simplex in $(X/G)^{\star T}$ the $N$ vertices of which are distinct $G$-orbits in $X$ mapped by $f \in \Sigma$ to some $\rho(G)$-orbits in $C\Lambda$. Taking a representative $x_j \in X$ of any of these $N$ \ $G$-orbits, and the corresponding points $f(x_j) \in C\Lambda$, we obtain a point of the direct product $X^N \times (C\Lambda)^N$. The orbit of this point under the action of all elements (\ref{semi}) described by (\ref{semact}) is a point of the space $\Pi_{N,G}(X)$ that does not depend on the choice of the representatives $x_j$.

By the construction and the first property of $k$-sufficient approximations, the resulting map $\sigma_N \setminus \sigma_{N-1} \to \Pi_{N,G}(X)$ is a fiber bundle, whose fiber is the set of all points $(s,f)$ such that $s$ is in some open simplex with $N$ vertices in $(X/G)^{\star T}$, and $f$ belongs to the affine subspace in ${\mathcal F}$ that consists of the maps which take fixed values on a fixed set of $N$ different $G$-orbits in $X$. The orientation rules for these bundles of simplices and subspaces follow immediately from the construction. \hfill $\Box$ 

\begin{corollary} \label{cormain}
If ${\mathcal F}$ is a $k$-sufficient subspace, then there is a spectral sequence $E^{p,q}_r({\mathcal F})$ converging to the group $H^*({\mathcal F} \setminus \Sigma)$ with the term $E^{p,q}_1$ given by formula $($\ref{absvf}$)$ for $p \geq -k$ and equal to zero if either $p>0$ or $q < -p(W-\dim X + \dim G - \dim C\Lambda)$. In particular, formula $($\ref{absvf}$)$ gives one $E^{p,q}_1({\mathcal F})$ for all $p,q$ such that $p<0$ and $p+q < (k+1)(W-\dim X + \dim G - \dim C\Lambda -1)$.
\end{corollary}

\noindent
{\it Proof}. Our filtration $\{\sigma_N\}$ defines a spectral sequence $E^r_{N,t}({\mathcal F})$ calculating the group $\bar H_*(\sigma)$; in particular $E^1_{N,t}({\mathcal F}) \simeq \bar H_{N+t}(\sigma_N \setminus \sigma_{N-1})$. Let us turn this homological spectral sequence into a cohomological one by setting 
\begin{equation}
\label{inver}
E^{p,q}_r({\mathcal F}) \equiv E^r_{-p, \dim {\mathcal F}-q-1}({\mathcal F}).\end{equation} The resulting spectral sequence converges to the group $\tilde H^*({\mathcal F}\setminus \Sigma)$, which is Alexander dual to the group $\bar H_*(\Sigma \cap {\mathcal F})$ calculated by the previous spectral sequence. 
Statement 2 of Proposition \ref{pro6} implies that all its non-trivial groups $E^{p,q}_r({\mathcal F}),$ $r\geq 1$, lie in the wedge (\ref{wedge}), in particular such groups with $p<-k$ are trivial if $p+q<(k+1)(W-\dim X + \dim G - \dim C\Lambda -1)$. Statement 3 of Proposition \ref{pro6} implies that all groups $E^{p,q}_1({\mathcal F})$ with $p \geq -k$ are as described by formula (\ref{absvf}). 

The construction in the relative case and the proof of formula (\ref{relaf}) are the same with obvious modifications. \hfill $\Box$

\begin{corollary}
\label{cor78}
The group $H_i(\mbox{\rm EM}_G(X,A; \R^W) \setminus \Sigma)$ is finitely generated for any $i$.
\end{corollary}

\noindent{\it Proof.} Choose a natural $k$ such that $i <(k+1)(W-\dim X + \dim G -\dim C\Lambda-1)$, then by Corollary \ref{cormain}  the number of independent generators of the group $H_i({\mathcal F} \setminus \Sigma)$  for any $k$-sufficient subspace ${\mathcal F}$ is effectively bounded by a finite number $\mu$. If the group $H_i(\mbox{\rm EM}_G(X,A; \R^W) \setminus \Sigma)$ is not finitely generated, let us choose some $\mu+1$ of its independent generators. By Proposition \ref{appro} they can be realized by cycles in some finite-dimensional subspace, and hence also by cycles in a $k$-sufficient subspace approximating it: a contradiction. \hfill $\Box$

\subsection{Isomorphism and stabilization of spectral sequences}

In this section we prove that not only the groups $E^{p,q}_1({\mathcal F})$ for all $k$-sufficient subspaces ${\mathcal F}$ are isomorphic to one another in the stable domain of the $(p,q)$-plane, but the entire spectral sequences are isomorphic there, too.

\begin{proposition}
\label{cor17}
If $L\geq  k(W+\dim X - \dim G + \dim C \Lambda) +1   $ then
all spectral sequences $E^{p,q}_r({\mathcal F})$ corresponding to $k$-sufficient $L$-dimensional affine subspaces ${\mathcal F}$ of $ \mbox{\rm EM}_G(X,A; \R^W)$ are isomorphic to one another in the domain of the $(p,q)$-plane where $p+q < (k+1)(W-\dim X + \dim G - \dim C\Lambda -1)$. \hfill $\Box$ 
\end{proposition}

\noindent
{\it Proof.} The $k$-th term $\sigma_k({\mathcal F})$ of the simplicial resolution of $\Sigma \cap {\mathcal F}$ can be constructed as above for an arbitrary $k$-transversal (not necessarily $k$-sufficient) affine subspace ${\mathcal F}$, and its terms $\sigma_N({\mathcal F}) \setminus \sigma_{N-1}({\mathcal F})$ are described by statement 3 of Proposition \ref{pro6} as well. Denote by $E^r_{p,q}({\mathcal F},k)$ the homological spectral sequence approximating the group $\bar H_*(\sigma_k({\mathcal F}))$ and  defined by our standard filtration of this simplicial resolution.

\begin{lemma}
\label{fam}
Let $\{{\mathcal F}_\tau\}$, $\tau \in [0,1]$, be a continuous family of $L$-dimen\-si\-onal $k$-transversal subspaces in $ \mbox{\rm EM}_G(X,A; \R^W)$. Then all spectral sequences $E_{p,q}^r({\mathcal F}_\tau,k)$ are isomorphic to one another.
\end{lemma}

\noindent
{\it Proof of Lemma \ref{fam}}. Denote by $\hat \sigma_k$ the space of all pairs $(\tau,z)$ where $\tau \in [0,1]$ and $z$ is a point of the term $\sigma_k({\mathcal F}_\tau)$ of the simplicial resolution of $\Sigma \cap {\mathcal F}_\tau$. This space admits an obvious projection $\pi:\hat \sigma_k \to [0,1]$ defined by the first elements of these pairs, and is filtered by the subspaces $\hat \sigma_N$, $N =1, \dots, k$. Consider the homological spectral sequence defined by this filtration and converging to the group $\bar H_*(\hat \sigma_k)$. For any $\tau \in [0,1]$, the  inclusion $\sigma_k({\mathcal F}_\tau) \equiv \pi^{-1}(\tau) \cap \hat \sigma_k \hookrightarrow \hat \sigma_k$ defines a homomorphism of spectral sequences. By the construction, for any $N=1, \dots, k$ the space $\hat \sigma_N \setminus \hat \sigma_{N-1}$ is a locally trivial fiber bundle over the segment $[0,1]$ with fibers $\sigma_N({\mathcal F}_\tau) \setminus \sigma_{N-1}({\mathcal F}_\tau)$. Therefore this homomorphism of spectral sequences defines an isomorphism of their terms $E^1$, and hence it also is an isomorphism. In particular, the spectral sequences $E_{p,q}^r({\mathcal F}_\tau,k)$ associated with all values $\tau \in [0,1]$ are isomorphic to one another. \hfill $\Box$ \medskip

By statement 1 of Proposition \ref{ksuff}, in the conditions of Proposition \ref{cor17} the space of $k$-transversal $L$-dimensional affine subspaces of $ \mbox{EM} _G(X,A; \R^W)$
 is path-connected if $L$ is large enough, therefore the corresponding spectral sequences $E^r_{p,q}({\mathcal F},k)$ are all isomorphic to one another. However, for $k$-sufficient subspaces ${\mathcal F}$ these spectral sequences are isomorphic to the (non-restricted) spectral sequences $E^r_{p,q}({\mathcal F})$ in the domain where $p+q\geq L-(k+1)(W-\dim X + \dim G - \dim C\Lambda -1)$, hence these spectral sequences are isomorphic to one another too, and the corresponding cohomological spectral sequences $E^{p,q}_r({\mathcal F})$ (related with them by  (\ref{inver})) also are there isomorphic to one another in the domain where $p+q <  (k+1)(W-\dim X + \dim G - \dim C\Lambda -1)$. \hfill $\Box$ \medskip

Let us prove that these spectral sequences coincide  in such a domain also for $k$-sufficient subspaces ${\mathcal F}$ of different sufficiently large dimensions.

\begin{lemma}[cf. \cite{Book}]
\label{lem18}
Let ${\mathcal F} \subset {\mathcal F}'$ be two $k$-sufficient subspaces in $ \mbox{EM} _G(X,A; \R^W)$, and ${\mathcal F}$ be in general position with the semialgebraic variety $\Sigma \cap {\mathcal F}'$ $($that is, the closure of ${\mathcal F}$ in the projective compactification of ${\mathcal F}'$ is transversal to the stratified subvariety consisting of $\Sigma \cap {\mathcal F}'$ and the hyperplane added at the compactification$)$. Then the corresponding spectral sequences $E^{p,q}_r({\mathcal F})$ and $E^{p,q}_r({\mathcal F}')$ converging to the cohomology groups of ${\mathcal F} \setminus \Sigma$ and ${\mathcal F}' \setminus \Sigma$ are isomorphic to one another in the domain where $p+q < (k+1)(W-\dim X + \dim G - \dim C\Lambda -1)$, and the identical embedding ${\mathcal F} \setminus \Sigma \subset {\mathcal F}' \setminus \Sigma$ induces an isomorphism of these cohomology groups up to dimension $(k+1)(W-\dim X + \dim G - \dim C\Lambda -1)-1$.
\end{lemma}

\noindent
{\it Proof.}
By the Thom isotopy lemma (see e.g. \cite{GM}) there is a tubular neighborhood $U$ of ${\mathcal F}$ in ${\mathcal F}'$ such that the pair of stratified varieties $(U, \Sigma \cap U)$ is homeomorphic to the product of the pair $({\mathcal F}, \Sigma \cap {\mathcal F})$ and an open ball of dimension \ $\dim {\mathcal F}' -\dim {\mathcal F}$. The group $H^*(U \setminus \Sigma)$ can be calculated using a spectral sequence constructed in exactly the same way as $E^{p,q}_r({\mathcal F}')$, but with the terms $\sigma_N({\mathcal F}' ) \subset (X/G)^{\star T} \times {\mathcal F}'$ replaced by their intersections with the subset $(X/G)^{\star T} \times U$. The space $\sigma(U) $ of the simplicial resolution appearing in this construction is homeomorphic as a filtered space to the direct product of $\sigma({\mathcal F})$ and an open ball of dimension $\dim {\mathcal F}' -\dim {\mathcal F}$, hence the corresponding spectral sequence $E^{p,q}_r(U)$ is isomorphic to $E^{p,q}_r({\mathcal F})$. On the other hand, the inclusion $U \hookrightarrow {\mathcal F}'$ induces a natural homomorphism of the homological spectral sequences $E_{p,q}^r({\mathcal F}' ) \to E_{p,q}^r(U)$. It is an isomorphism of $E_{p,q}^1$ for $p \leq k$. Indeed, both terms $\sigma_p(U) \setminus \sigma_{p-1}(U)$ and $\sigma_p({\mathcal F}') \setminus \sigma_{p-1}({\mathcal F}')$ are the spaces of fiber bundles over one and the same base; the fibers of the latter bundle are some affine spaces, and fibers of the former one are open balls in these affine spaces.
All non-trivial terms $E_{p,q}^1$ of these spectral sequences with $p>k$ lie in the domain where $p+q \leq \dim {\mathcal F}'-(k+1)(W-\dim X + \dim G- \dim C\Lambda -1)-1$. Therefore neither these groups nor any higher differentials from them affect the groups $E_{p,q}^r$ with $p+q > \dim {\mathcal F}' - (k+1)(W-\dim X + \dim G - \dim C\Lambda -1)-1$ or, which is the same, the terms $E_r^{p,q}$ of the cohomological spectral sequence with $p+q < (k+1)(W-\dim X + \dim G - \dim C\Lambda -1).$ \hfill $\Box$ \medskip

Theorems \ref{absv} and \ref{relav} follow immediately from propositions \ref{appro}, \ref{ksuff}, \ref{pro6} (with Corollaries \ref{cormain} and \ref{cor78}), \ref{cor17}, and Lemma \ref{lem18}. Namely, for any natural number $T$ the desired spectral sequence $E^{p,q}_r$ coincides in the domain where $p+q < T$ with the spectral sequence  $E^{p,q}_r({\mathcal F})$ constructed in \S \ref{sre}, where ${\mathcal F}$ is an arbitrary $k$-sufficient affine subspace in $ \mbox{EM} _G(X,A; \R^W)$, the dimension of which and number $k$ are large enough with respect to $T$.

\subsection{Which equivariant function spaces can be reduced to the previous pattern?}
\label{which}

In the case of the trivial group $G$, the previous scheme allows one to calculate the cohomology groups of the spaces of continuous maps $X \to Y$ where $X$ is $m$-dimensional and $Y$ is an $m$-connected finite CW-complex. Indeed, we can embed $Y$ into a sphere $S^{W-1}$ as a deformation retract of the complement of a subcomplex $\Lambda$ of codimension $\geq m+2$ in $S^{W-1}$, so that the spaces of maps to $Y$ and to $\R^W \setminus C\Lambda$ are homotopy equivalent to one another. However, the proof of this fact known to me (see \cite{Book}) does not immediately extend to the case of non-trivial group actions. 

On the other hand, it seems likely that even if we cannot realize this strategy geometrically, it can sometimes give us a hint about how a spectral sequence calculating the desired cohomology groups could look like, and the limit cohomology group of this virtual spectral sequence could be a correct guess at the answer. More precisely, suppose that $X$ is $m$-dimensional, $Y$ is $m$-connected, $G$ is a finite group acting freely on both $X$ and $Y$, and we study the homology of the space of $G$-equivariant maps $X \to Y.$

It is easy to construct a representation
 $\rho: G \to \mbox{O}(\R^W)$ with sufficiently large $W$ and also an embedding of $Y$ into the unit sphere $S^{W-1} \subset \R^W$ commuting with the actions of $G$. For instance, we can take for $\R^W$ a sufficiently large Cartesian power of the space of the regular representation of $G$. Increasing this power, we can make the codimension of the union of non-regular orbits arbitrarily large. The complement of this union in $S^{W-1}$ together with the $G$-action on it can then be made to homotopically approximate arbitrarily well the universal covering of $K(G,1)$. If the codimension of this union of non-regular orbits is greater than \ $\dim Y$, then there are no obstructions to embedding the space $Y/G$ into the quotient of this complement by this $G$-action so that after lifting the embedding to the $G$-covering we get a $G$-equivariant embedding $Y \to S^{W-1}$.

To solve our problem using the spectral sequence described above, we need to find a compact subcomplex $\Lambda \subset S^{W-1} \setminus Y$ invariant under our $G$-action and such that

a) $\Lambda$ contains all non-regular orbits of our $G$-action;

b) the complex $Y \subset S^{W-1}$ is a $G$-equivariant deformation retract of $S^{W-1} \setminus \Lambda$;

c) the codimension of $\Lambda$ in $S^{W-1}$ is at least $m +1$. 

The main difficulty is with the property (c). The proof that for trivial $G$ there exists a $\Lambda$ with properties (a)-(c) (see 
 \cite{Book}) does not seem to generalize to $G \neq \{1\}$ in a straightforward way. Indeed, if we try to apply it, we need to cancel the critical points with high Morse indices of $G$-invariant functions on the manifold $S^{W-1} \setminus Y$. This is essentially equivalent to doing the same for usual functions on the space of regular $G$-orbits, which is not simply-connected if $G$ is not connected, and the standard methods (see e.g. \cite{hcob}) may not work. Whether or not such a representation $\rho: \mbox{O}(\R^W)$ and a subcomplex $\Lambda$ exist for arbitrary $G$ seems to be a good problem.

However, even if we are unable to find $\rho$ and $\Lambda$ with these properties, let us take $\Lambda$ to be an arbitrary subcomplex that only satisfies conditions (a) and (b), and then construct the complexes $\Pi_{N,G}(X)$ and write down our spectral sequence exactly as in the previous subsections. By the Alexander duality, the {\it homological} dimension of $\Lambda$ does not exceed $W-m-2$, therefore the term $E^{p,q}_1$ of our spectral sequence looks as previously, in particular it has only finitely many non-trivial groups on any diagonal of the form $\{p+q= \mbox{const} \}$, and all these groups are finitely generated. It seems likely that this spectral sequence still will calculate the desired cohomology groups of the space of equivariant maps, even though our geometrical realization of it does not work in this case.

\medskip
I thank F.~Cohen and V.~Turchin for advices and conversations.

\end{document}